\documentclass[11pt]{article}

\usepackage{amsthm}

\newcommand\tsum{\textstyle\sum\nolimits}
\newcommand{\abs}[1]{\left|#1\right|}

\usepackage[nocompress]{cite}
\usepackage{graphicx}

\usepackage{hyperref}
\hypersetup{
    colorlinks=true,
    linkcolor=blue,
    anchorcolor=blue,
    citecolor=blue,
}

\usepackage{enumitem}

 \usepackage{amsmath}
\usepackage{amsfonts}
\usepackage{amssymb, amsfonts, amsmath, nicefrac}
\usepackage{color}

\usepackage[margin=1in]{geometry}

\newcommand{\cQ}{{\mathfrak Q}}

\newcommand{\cP}{{\mathfrak P}}
\newcommand{\cD}{{\mathfrak D}}
\newcommand{\cM}{{\mathfrak M}}

\newcommand{\cC}{{\mathfrak C}}

\newcommand{\cR}{{\mathfrak R}}

\newcommand{\ch}{{\mathfrak h}}

\newtheorem{theorem}{Theorem}[section]
\newtheorem{lemma}{Lemma}[section]
\newtheorem{proposition}{Proposition}[section]
\newtheorem{example}{Example}[section]
\newtheorem{remark}{Remark}[section]
\newtheorem{definition}{Definition}[section]

\newtheorem{assumption}{Assumption}[section]
\newtheorem{cor}{Corollary}[section]

\def\argmin{\mathop{\rm arg\,min}}
\def\argmax{\mathop{\rm arg\,max}}

\newcommand{\Q}{{\cal Q}}

\newcommand{\V}{{\cal V}}
\newcommand{\Z}{{\cal Z}}
\newcommand{\F}{{\cal F}}

\newcommand{\C}{{\cal C}}
\newcommand{\A}{{\cal A}}

\newcommand{\X}{{\cal X}}

\newcommand{\PP}{{\cal P}}

\newcommand{\Y}{{\cal Y}}

\newcommand{\W}{\mathcal{W}}
\newcommand{\s}{{\cal S}}
\newcommand{\R}{{\cal R}}

\newcommand{\be}{\begin{equation}}
\newcommand{\ee}{\end{equation}}

\def\w{\omega}

\def\O{\Omega}

\newcommand{\avr}{{\sf AV@R}}

\def\bbr{{\Bbb{R}}} 
\def\bbe{{\Bbb{E}}} 

\def\bbp{{\Bbb{P}}}

\begin{document}

\begin{titlepage}
\title{\bf Rectangularity and duality of distributionally robust Markov Decision Processes}
\author{{\bf Yan Li}\thanks{
Georgia Institute of Technology, Atlanta, Georgia
30332, USA, \tt{yli939@gatech.edu}} \and  {\bf Alexander Shapiro}\thanks{Georgia Institute of Technology, Atlanta, Georgia
30332, USA, \tt{ashapiro@isye.gatech.edu}\newline
Research of this
author was partially supported by Air Force Office of Scientific Research (AFOSR)
under Grant FA9550-22-1-0244.}
}


\maketitle
\begin{abstract}
The main goal of this paper is to discuss several approaches to the formulation of  distributionally robust counterparts of Markov Decision Processes,  where the transition kernels are not specified exactly but rather are assumed  to be elements of the corresponding ambiguity sets.
The intent is to clarify some connections between the game and static formulations of distributionally robust MDPs,  and delineate the role of rectangularity associated with ambiguity sets in determining these connections.
\end{abstract}

{\bf Keywords:} Markov Decision Process, distributional robustness, game formulation,  strong duality, risk measures

\end{titlepage}

\setcounter{equation}{0}
\section{Introduction}
\label{sec-intr}
Consider a (finite horizon) Markov Decision Process (MDP)
\begin{equation}\label{intra-1}
 \min_{\pi\in \Pi} \bbe^{\pi}
 \left [ \tsum_{t=1}^{T}
c_t(s_t,a_t,s_{t+1})+c_{T+1}(s_{T+1})
\right].
\end{equation}
Here $\s_t$ is the state space,  $\A_t$ is the action set and $c_t:\s_t\times \A_t\times \s_{t+1}\to \bbr$ is the cost function,
at stage $t=1,...,T$. For the sake of simplicity, in order to concentrate on the main issues discussed in this paper, we assume that the sets $\s_t$ and $\A_t$ are finite. Unless stated otherwise the optimization in \eqref{intra-1} is performed over randomized   Markovian  policies $\pi_t:\s_t\to \Delta_{|\A_t|}$, $t=1,...,T$, where $|\A_t|$ is the cardinality of the set  $\A_t$ and $\Delta_n=\{q\in \bbr^n_+:\tsum_{i=1}^n q_i=1\}$   denotes the $n$-dimensional simplex.
We denote the set of randomized Markovian policies as $\Pi$.
We also allow the action set $\A_t(s_t)$ to depend on $s_t\in \s_t$, $t=1,...,T$, and assume the set $\A_t(s_t)$ is {\em nonempty} for every $s_t\in \s_t$,
and that  $s_1 \in \s_1$ is a fixed state.
It is also possible to consider random $s_1$ following a given distribution without changing the essence of our discussions.
For a policy $\pi\in \Pi$,  the expectation
$\bbe^{\pi}$ is taken with respect to the probability law defined by the transition kernels $\{P_t\}$.
 Since both the state and action spaces are finite,
 we can view a kernel $P_t$ as a finite dimensional vector with components $P_t(s_{t+1}|a_t,s_t)$, $s_{t+1}\in \s_{t+1}$, $a_t\in \A_t(s_t)$, $s_t\in \s_t$.
In particular, $P_t(\cdot|s_t, a_t)$ specifies the probability distribution of $s_{t+1}$ conditioned on the current state-action pair $(s_t, a_t)$.

The main goal of this paper is to discuss approaches to the formulation of  distributionally robust counterparts of the MDP \eqref{intra-1},  where the transition kernels $\{P_t\}$ are not specified exactly but rather are assumed  to be elements of the corresponding ambiguity sets.
Existing literature largely focuses on the so-called static formulation\footnote{What is called in this article  the static formulation, in some publications   is  referred to as {\em robust } MDPs, e.g., \cite{iyen,Nilim2005,XuMannor2012}.}, where the objective in \eqref{intra-1} is replaced by the worst-case scenario when the transition kernel defining the expectation is subject to variations.
The tractability of the static formulation hinges upon the existence of certain robust version of dynamic equations. To certify such dynamic equations,  structural conditions on the ambiguity sets, termed rectangularity, have been identified in prior studies \cite{goyal2023robust, goh2018data, iyen,Nilim2005, Tallec, kuhn2013}.
 Oftentimes, establishing dynamic equations with different classes of ambiguity sets  follows a case-by-case analysis, 
 e.g., for $\mathrm{(s,a)}$-rectangularity \cite{iyen,Nilim2005}, $\mathrm{s}$-rectangularity \cite{kuhn2013}, and $\mathrm{r}$-rectangularity \cite{goyal2023robust, goh2018data}.

A game formulation is identified for $(\mathrm{s,a})$-rectangular sets \cite{iyen,Nilim2005}, where a turn-based zero-sum game with perfect information is formulated.
It is shown therein that dynamic equations naturally arise once the sequential game is constructed.  
Compared to its static counterpart, game formulations of  distributionally robust MDPs have received limited attention.
In addition, the aforementioned turn-based zero-sum game is asymmetric as the information of two players is different.
The construction of this asymmetric game also appears to rely on the $(\mathrm{s,a})$-rectangularity of the ambiguity set.
Asymmetric games with a class of non-rectangular sets have also been considered  in \cite{mannor2012lightning, mannor2016robust}, which we discuss in Section \ref{subsec_static_form}.

In this manuscript, we consider a natural game formulation of distributionally robust MDPs.
This indeed serves as a unifying approach to both the game and the static formulations under a variety of rectangular structures of the ambiguity set.
We establish dynamic equations and strong duality for the game formulation. Notably these results do not necessarily require rectangularity of the ambiguity sets.
Then equivalence between game and static formulations is established, under a condition that (through simple verification) covers all existing mainstream rectangular sets (e.g., \cite{goyal2023robust, goh2018data, iyen,Nilim2005, Tallec, kuhn2013}).
Through this equivalence, dynamic equations and strong duality for the static formulation immediately follows.
We discuss further implications of this equivalence and propose a class of rectangular sets that encompasses existing ones.
Connections to risk-averse setting are subsequently made, and we highlight an essential difference when considering the Optimal Control framework in the distributionally robust setting.

 This paper is organized as follow.
In Section \ref{sec-boc}, we introduce the game formulation and its dual,  and establish their dynamic equations through an elementary argument.
From dynamic equations we naturally investigate the strong duality of the game formulation.
We propose two structural assumptions on the ambiguity sets that serve as sufficient conditions of the strong duality, and discuss their implications in  determining the existence of non-randomized optimal policies.
We discuss how these assumptions can be satisfied by existing rectangularity sets.
On the flip side, we show that these structural assumptions do not always require rectangularity of the ambiguity sets.
We conclude this section by explicitly discussing history-dependent policies in Subsection \ref{histdepen}.

In Section \ref{subsec_static_form}, we show the dynamic equations of the game formulation provide a useful tool to establish equivalence between the game and static formulations, thus facilitating the study of the latter.
As immediate consequences of this equivalence: (i)
 dynamic equations carry to the static formulations without additional overhead;
(ii)  strong duality of the static formulation is also established.
Notably, this allows one to consider static formulations with different ambiguity sets with a unifying treatment.
As a byproduct, it is shown that any distributionally robust MDP with an $\mathrm{r}$-rectangular set can be equivalently reformulated into another one with a $\mathrm{s}$-rectangular set. We provide an elementary construction of such a reformulation. 
Towards the end of the section, we construct a class of ambiguity set termed as $\mathrm{sr}$-rectangular sets, and derive dynamic equations and strong duality for the static formulation as an application of our analytical approach.
The proposed rectangularity contains strictly the  $\mathrm{sa}$- \cite{iyen,Nilim2005}, $\mathrm{s}$- \cite{Tallec, kuhn2013}, and $\mathrm{r}$-rectangular \cite{goyal2023robust, goh2018data} sets. 

In Section \ref{sec_cost_robust}, we show all the aforementioned results hold similarly for cost-robust MDP, where cost functions, instead of the transition kernels, are subject to ambiguity.
Notably, the discussion therein follows exactly the same argument as for the distributionally robust MDP.
As a byproduct, a cost-robust MDP can be reduced to a regular (i.e., non-robust) MDP with a convex policy-dependent cost function.

In Section \ref{sec_risk_averse}, we consider risk-averse setting, where we discuss a natural construction of ambiguity sets based on law invariant coherent risk measures.
Finally,  in Section \ref{sec-soc}, we point out to an essential difference between distributionally robust formulations in the MDP and Optimal Control  frameworks.

Note that
since the set $\s_t$ is finite, a probability distribution supported on
$\s_t$ can be identified with probability vector $q\in \Delta_{|\s_t|}$.
By $\delta_a$ we denote  Dirac measure  of mass one at  point $a$. 
Let us recall the following properties of the   min-max problem:
 \begin{equation}\label{minmax-1}
\min_{x\in \X}\sup_{y\in \Y} f(x,y),
\end{equation}
 where $\X$ and $\Y$ are nonempty sets and $f:\X\times \Y\to \bbr$ is a real valued function.
A point  $(x^*,y^*)\in \X\times \Y$ is a saddle point of   problem \eqref{minmax-1}
if
$
x^*\in \argmin_{x\in \X} f(x,y^*)$ and $y^*\in \argmax f(x^*,y).
$
If the saddle point exists, then problem \eqref{minmax-1} has the same optimal value as its dual
\begin{equation}\label{minmax-2}
\max_{y\in \Y} \inf_{x\in \X} f(x,y),
\end{equation}
and $x^*$ is an optimal solution of problem \eqref{minmax-1} and $y^*$ is an optimal solution of problem \eqref{minmax-2}. Conversely,
if the    optimal values of problems  \eqref{minmax-1} and  \eqref{minmax-2} are the same, and $x^*$ is an optimal solution of problem \eqref{minmax-1} and $y^*$ is an optimal solution of problem \eqref{minmax-2}, then
$(x^*,y^*)$ is a saddle point.
If $\X$ and $\Y$ are convex subsets of finite dimensional vector spaces, $f(x,y)$ is continuous, convex in $x$ and concave in $y$, and at least one of  the sets $\X$ or  $\Y$ is compact, then the optimal values of problems  \eqref{minmax-1} and  \eqref{minmax-2} are equal to each other. This is a particular case of  Sion's minimax theorem \cite{sion}.

\setcounter{equation}{0}
\section{Distributionally robust  MDP}
\label{sec-boc}

Assume that for $t=1,...,T$, there is a set $\PP_t$ of transition kernels $P_t$.
At this point we do not  specify a particular construction of the ambiguity  sets $\PP_t$. Unless stated otherwise we make the following assumption:
{\em the sets $\PP_t$,  are  closed}.
Since the   probabilities $P_t(s_{t+1}|a_t,s_t)\in [0,1]$, the sets $\PP_t$, $t=1,...,T$, are bounded and hence are compact.

\subsection{Game formulation of distributionally robust MDP}\label{subsec_game_form}

The game formulation of distributionally robust MDPs considers  a dynamic game between the decision maker (the controller) and the  adversary (the nature).
The nature can choose a kernel $P_t \in \PP_t$ for each $s_t \in \s_t$, $t=1,...,T$,
 and this defines a   policy of the nature (cf., \cite{Jask2016, Shapley}).
 Consider the decision process determined by the    decision  history
 $s_1,a_1,P_1,s_2,...,a_{t-1},P_{t-1}, s_t$, $t=1,...,T,$   where $s_t\in \s_t$,  $a_t\in \A_\tau(s_t)$, $P_t\in \PP_t$.
 The controller  chooses a randomized   policy of the form $\pi_{t}(\cdot|s_{t},\ch_{t-1}): \s_t \times \mathfrak{H}_{t-1} \to  \Delta_{|\A_t(s_{t})|}$,    $t=1,...,T$, where
\begin{equation}\label{history}
 \ch_{t-1}:=(s_1,a_1,P_1,...,s_{t-1},a_{t-1},P_{t-1}),
\end{equation}
and $\mathfrak{H}_{t-1}$ denotes the set of possible $\ch_{t-1}$.
The nature chooses its policy $\gamma_t(s_t,\ch_{t-1}): \s_t \times \mathfrak{H}_{t-1} \to \PP_t$, $t = 1, \ldots, T$.
  It is said that  a  policy $\{\pi_t\}$ is non-randomized if for every $s_t \in \s_t$, $t=1,...,T$,   the corresponding  probability distribution  $\pi_{t}(\cdot|s_{t},\ch_{t-1})$ is supported on a single point of $\A_t(s_t)$  (i.e., is the delta function).
 Note that for a given pair of policies $(\pi_t, \gamma_t)$, the controller and the nature choose their respective actions simultaneously at every $s_t \in \s_t$.

 Policies $\{\pi_t\}$ and $\{\gamma_t\}$, of the controller and nature, define the respective probability distribution  on the set of the histories of the decision process (Ionescu Tulcea theorem). The corresponding expectation is denoted $\bbe^{\pi,\gamma}$.
 It is said that a policy $\{\pi_t\}$ of the controller is  Markovian  if
$\pi_{t}(\cdot|s_{t})$,    $t=1,...,T$, does not depend on the history $\ch_{t-1}$. Similarly a policy of the nature is Markovian if
$\gamma_t(s_t)$ is a function of $s_t$ alone. Unless stated otherwise we deal with Markovian policies of the controller and the  nature\footnote{History-dependent policies are discussed in Section  \ref{histdepen}.}, and denote by
$\Pi$  the set of (randomized) {\em Markovian} policies  of the controller and by $\Gamma $  the set of  {\em  Markovian} policies  of the nature.
The corresponding  problem of the game formulation is

\begin{equation}\label{mdp-1}
 \min_{\pi\in \Pi}\max_{\gamma\in \Gamma}\bbe^{\pi,\gamma}
 \left [ \tsum_{t=1}^{T}
c_t(s_t,a_t,s_{t+1})+c_{T+1}(s_{T+1})
\right].
\end{equation}
We refer to \eqref{mdp-1} as the {\em primal} problem of the game formulation.
In the primal problem, the nature chooses its policy after the controller.

 It is well known that in the case of the   MDP problem \eqref{intra-1},  the optimal policies are Markovian (cf., \cite{puterman2014markov}).
We will discuss in Section \ref{histdepen}
  the issue  of non-Markovian optimal policies in the distributionally robust setting.

Given a policy $\{\pi_t\}$ of the controller, the nature chooses a policy
so as to maximize the total expected cost, i.e.,
\begin{equation}\label{mdp-2}
\max_{\gamma \in \Gamma}  \bbe^{\pi, \gamma}\left [ \tsum_{t=1}^{T}
c_t(s_t,a_t,s_{t+1})+c_{T+1}(s_{T+1})
\right].
\end{equation}

\begin{proposition}
\label{pr-dyncontr}
Given a policy $\pi\in \Pi$ of the controller, problem
 \eqref{mdp-2} admits the following
 dynamic programming equations: $\V^\pi_{T+1}(s_{T+1})=  c_{T+1}(s_{T+1})$ and for $t=T,...,1$ and $s_t\in \s_t$,
 \begin{align}\label{mdp-3}
\V^\pi_t(s_t)  = \max_{P_t \in \PP_t}  \underbrace{\tsum_{s_{t+1} \in \s_{t+1}}  \tsum_{a_t \in \A_t(s_t)} P_t(s_{t+1} | s_t ,a_t) \pi_t(a_t | s_t) \big[
c_t(s_t, a_t ,s_{t+1})+\V^\pi_{t+1}(s_{t+1}) \big]}_{ \bbe^{\pi_t, P_t}[
c_t(s_t, a_t ,s_{t+1})+\V^\pi_{t+1}(s_{t+1})]  },
\end{align}
where   $\V^\pi_t(s_t)$  equals the optimal value of problem
\begin{align}\label{eq_inner_max_starting_t}
\max_{\gamma \in \Gamma}  \bbe^{\pi, \gamma}\left [ \tsum_{i=t}^{T}
c_i(s_i,a_i,s_{i+1})+c_{T+1}(s_{T+1})  \big| s_t
\right].
\end{align}
\end{proposition}

\begin{remark}
{\rm
  By writing $P_t\in \PP_t$  we mean that the transition kernel $P_t$ is viewed as an {\em element} of the ambiguity set $\PP_t$.
  Throughout the rest of our discussion,  a particular choice of $P^*_t\in \PP_t$ (for instance, the optimal solution of \eqref{mdp-3}) can depend on the state $s_t\in \s_t$, i.e., $P^*_t$ can be different for different values of $s_t$.
We will say explicitly that $P^*_t$ is independent of $s_t$ when the considered element  $P^*_t$ of the ambiguity set $\PP_t$    is  the same  for all $s_t\in \s_t$. Of course the transition probability  $P_t^*(\cdot|s_t,a_t)$
 depends on $s_t$.
 }
 $\hfill \square$
\end{remark}

\begin{proof}
Solving \eqref{mdp-2} is equivalent to solving a regular MDP of the nature defined as follows.
At stage $t$, the state space is given by $\s_t$;
the possible actions at $s_t \in \s_t$ are given by $\PP_t$;
the cost function of the nature is given by $-c_t$; and the transition probability $\mathbb{P}^{\mathrm{nature}}_t$ is given by
\begin{equation}\label{tranat}
\mathbb{P}^{\mathrm{nature}}_t (s_{t+1} |s_t , P_t) = \tsum_{a_t \in \A_t(s_t)} P_t(s_{t+1} | s_t, a_t) \pi_t(a_t | s_t).
\end{equation}
Note that the action space of nature's MDP is compact.
Applying standard dynamic equations (cf., Theorem 4.3.2, \cite{puterman2014markov}) of regular MDPs yields \eqref{mdp-3}.
\end{proof}

We also have that for a   policy  $\pi\in \Pi$, the corresponding optimal policy $\gamma^\pi_t(s_t)=P_t^\pi$
of the nature is determined  by dynamic equations
\begin{align}\label{mdp-3opt}
P_t^\pi \in \argmax_{P_t \in \PP_t}  \bbe^{\pi_t, P_t}[
c_t(s_t, a_t ,s_{t+1})+\V^\pi_{t+1}(s_{t+1})].
\end{align}
 The maximum of \eqref{mdp-2} is then minimized over $\pi\in \Pi$ and
 we can write  the following  dynamic programming equations for the primal problem \eqref{mdp-1}.

\begin{proposition}\label{prop_game_primal_dp}
The dynamic programming equations for the primal  problem \eqref{mdp-1} are: $V_{T+1}(s_{T+1})=  c_{T+1}(s_{T+1})$ and for $t=T,...,1$ and $s_t\in \s_t$,
\begin{equation}\label{mdp-5}
V_t(s_t)=\min_{\pi_t(\cdot|s_t) \in \Delta_{|\A_t(s_t)|}} \max_{P_t \in \PP_t} \tsum\limits_{s_{t+1} \in \s_{t+1}} \tsum\limits_{a_t \in \A_t(s_t)} P_t(s_{t+1} | s_t ,a_t) \pi_t(a_t | s_t) \big[
c_t(s_t, a_t ,s_{t+1})+V_{t+1}(s_{t+1}) \big],
\end{equation}
where  $V_t(s_t)$  equals the optimal value of
 \begin{align}\label{eq_primal_starting_at_t}
 \min_{\pi \in \Pi}
 \max_{\gamma \in \Gamma}  \bbe^{\pi, \gamma}\left [ \tsum_{i=t}^{T}
c_i(s_i,a_i,s_{i+1})+c_{T+1}(s_{T+1})  \big| s_t
\right].
 \end{align}
\end{proposition}

 \begin{proof}
 Fix a  policy $\pi \in \Pi$.
 We proceed to show that $V_t (s_t)\leq \V^\pi_t(s_t)$ by induction  for $t = T+1, \ldots 1$,
 where $\V^\pi_t(s_t)$ is the solution of \eqref{mdp-3}.
 The claim holds trivially at $t = T+1$.
 For any $t \leq T$, it holds that
 \begin{align*}
 V_t (s_t) & \overset{(a)}{\leq}
 \min_{\pi'_t(\cdot|s_t) \in \Delta_{|\A_t(s_t)|}} \max_{P_t \in \PP_t} \tsum\limits_{s_{t+1} \in \s_{t+1}} \tsum\limits_{a_t \in \A_t(s_t)} P_t(s_{t+1} | s_t ,a_t) \pi'_t(a_t | s_t) \big[
c_t(s_t, a_t ,s_{t+1})+ \V^\pi_{t+1}(s_{t+1}) \big] \\
& \leq
 \max_{P_t \in \PP_t} \tsum\limits_{s_{t+1} \in \s_{t+1}} \tsum\limits_{a_t \in \A_t(s_t)} P_t(s_{t+1} | s_t ,a_t) \pi_t(a_t | s_t) \big[
c_t(s_t, a_t ,s_{t+1})+ \V^\pi_{t+1}(s_{t+1}) \big]  \\
& \overset{(b)}{=} \V^{\pi}_t(s_t),
 \end{align*}
 from which we complete the induction step.
 Here $(a)$ follows from \eqref{mdp-5} and the induction hypothesis that $V_{t+1} (s_{t+1})\leq \V^\pi_{t+1} (s_{t+1})$,
 and $(b)$ follows from \eqref{mdp-3}.

 Let $V^*_t(s_t)$ denote the optimal value of \eqref{eq_primal_starting_at_t}.
 From \eqref{eq_inner_max_starting_t}, we have
 $V^*_t (s_t) = \min_{\pi \in \Pi} \V^\pi_t(s_t)$.
 Given that $V_t (s_t)\leq \V^\pi_t(s_t)$, taking minimum over $\pi \in \Pi$ then yields
 $V_t (s_t)\leq V^*_t(s_t)$.
 Finally, let $\pi^*$ be defined by \eqref{mdp-5},
 then from Proposition \ref{pr-dyncontr} it clearly holds that
 $V_t(s_t) = \V^{\pi^*}_t (s_t)$ and subsequently $V_t (s_t)\geq V_t^*(s_t)$.
Hence we obtain
 $V_t (s_t)= V^*_t(s_t)$, this   concludes  the proof.
 \end{proof}

Let us consider the   (Markovian)   policy
\begin{equation}\label{mark}
\pi^*_t(\cdot|s_t)\in  \argmin_{\pi_t(\cdot|s_t) \in \Delta_{|\A_t(s_t)|}} \max_{P_t \in \PP_t} \tsum\limits_{s_{t+1} \in \s_{t+1}} \tsum\limits_{a_t \in \A_t(s_t)} P_t(s_{t+1} | s_t ,a_t) \pi_t(a_t | s_t) \big[
c_t(s_t, a_t ,s_{t+1})+V_{t+1}(s_{t+1}) \big],
\end{equation}
$t=1,...,T$, determined by the dynamic programming equations \eqref{mdp-5}.
Since $\Delta_{|\A(s_t)|}$ is compact, the minimization problem in the right hand side of \eqref{mark} always has an optimal solution $\pi_t^*$. In view of Proposition \ref{prop_game_primal_dp},
the corresponding policy $\{\pi_t^*\}$ is an optimal solution (optimal policy of the controller) of the primal problem \eqref{mdp-1}. The non-randomized policies of the controller  are determined by the counterpart of    \eqref{mark} with $\pi^*_t(a_t | s_t) =\delta_{a^*_t(s_t)}$   being the Dirac measure, that is
\begin{equation}\label{mark-2}
a^*_t(s_t)\in  \argmin_{a_t \in  \A_t(s_t)} \max_{P_t \in \PP_t} \tsum\limits_{s_{t+1} \in \s_{t+1}} P_t(s_{t+1} | s_t ,a_t)\big[
c_t(s_t, a_t ,s_{t+1})+V_{t+1}(s_{t+1}) \big].
\end{equation}

The dual of the min-max problem in the right hand side of \eqref{mdp-5} is obtained by the interchange of the min and max operations. That is, given a policy $\gamma \in \Gamma$ of the nature, the corresponding dynamic programming equations for the controller are: $\Q^\gamma_{T+1}(s_{T+1})=  c_{T+1}(s_{T+1})$ and for $t=T,...,1$,
\begin{equation}\label{dual-3}
\Q^\gamma_t(s_t)=\min_{\pi_t(\cdot|s_t) \in \Delta_{|\A_t(s_t)|}}  \tsum\limits_{s_{t+1} \in \s_{t+1}} \tsum\limits_{a_t \in \A_t(s_t)} P_t(s_{t+1} | s_t ,a_t) \pi_t(a_t | s_t) \big[
c_t(s_t, a_t ,s_{t+1})+\Q^\gamma_{t+1}(s_{t+1}) \big].
\end{equation}
Consequently
the  dynamic programming equations for the dual  problem are: $Q_{T+1}(s_{T+1})=  c_{T+1}(s_{T+1})$ and for $t=T,...,1$,
  \begin{equation}\label{dual-4}
Q_t(s_t)=\max_{P_t\in \PP_t} \min_{\pi_t(\cdot|s_t) \in \Delta_{|\A_t(s_t)|}}   \tsum\limits_{s_{t+1} \in \s_{t+1}} \tsum\limits_{a_t \in \A_t(s_t)} P_t(s_{t+1} | s_t ,a_t) \pi_t(a_t | s_t) \big[
c_t(s_t, a_t ,s_{t+1})+Q_{t+1}(s_{t+1}) \big].
\end{equation}
The dynamic equation \eqref{dual-4}  define the dual
of the primal problem \eqref{mdp-1}. That is, starting with $Q_{T+1}(s_{T+1})=  c_{T+1}(s_{T+1})$ and going backward in time, eventually $Q_1(s_1)$ gives the optimal value of the dual problem. The optimal solutions $P_t^*\in \PP_t$, $t=1,...,T$,  in the right hand side of \eqref{dual-4} determine an optimal (Markovian) policy of the nature in the dual problem. Since the sets $\PP_t$ are assumed to be closed and hence compact, the maximum in \eqref{dual-4} is attained.

We write the dual problem as
\begin{equation}\label{dual-1}
 \max_{\gamma\in \Gamma} \min_{\pi\in \Pi} \bbe^{\pi,\gamma}
 \left [ \tsum_{t=1}^{T}
c_t(s_t,a_t,s_{t+1})+c_{T+1}(s_{T+1})
\right].
\end{equation}
In the dual setting the nature first chooses its policy $\gamma \in \Gamma$. Then  the controller finds the corresponding optimal (Markovian) policy $\pi \in \Pi$ to minimize the total expected cost, i.e.,
\begin{equation}\label{dual-2}
\min_{\pi\in \Pi} \bbe^{\pi, \gamma}
 \left [ \tsum_{t=1}^{T}
c_t(s_t,a_t,s_{t+1})+c_{T+1}(s_{T+1})
\right].
\end{equation}

By the standard theory of min-max problems we have that the optimal value of the dual problem \eqref{dual-1}  is less than or equal to the optimal value of the primal problem \eqref{mdp-1}, that is
\begin{equation}\label{opt-1}
 \mathrm{OPT}(\ref{dual-1})\le  \mathrm{OPT}(\ref{mdp-1}).
\end{equation}
We address now the question when the equality in \eqref{opt-1} holds, i.e., there is no duality gap between the primal and dual problems.
 In particular, the following remark illustrates our essential strategy.

\begin{remark}\label{remark_strong_dual_general}
{\rm
For $t=T+1$ we have that  $V_{T+1}(\cdot)=Q_{T+1}(\cdot)=c_{T+1}(\cdot)$. Going backward in time suppose that $V_{t+1}(\cdot)=Q_{t+1}(\cdot)$ and consider the dynamic programming equations \eqref{mdp-5}  and \eqref{dual-4} for the stage $t$. For  $s_t\in \s_t$ the function in the right hand side of \eqref{mdp-5}  can be considered as a function of two variables $P_t \in \PP_t$  and $\pi_t(\cdot|s_t)\in \Delta_{|\A_t(s_t)|}$.
Then we can talk about the concept of saddle point for the min-max problem \eqref{mdp-5} and its max-min dual problem \eqref{dual-4}.
Again, from weak duality we have
  \begin{equation}\label{dual-5}
Q_t(s_t)\le V_t(s_t),\;\;s_t\in \s_t,\;t=1,...,T.
\end{equation}
 If the saddle point exists for \eqref{mdp-5}, then based on  \eqref{dual-4} and our induction hypothesis that  $V_{t+1}(\cdot)=Q_{t+1}(\cdot)$, we have equality in \eqref{dual-5}.
Therefore if we can ensure existence of the saddle point for all $s_t\in \s_t$  and $t=1,...,T$,
then  going backward in time  eventually $V_1(\cdot)=Q_1(\cdot)$, and hence the optimal values of the primal problem \eqref{mdp-1} and its dual \eqref{dual-1} are equal to each other, i.e. the strong duality holds.
} $\hfill \square$
\end{remark}

Going forward we   introduce  conditions (Assumptions \ref{assum} and \ref{assump_margin_convex}) that   certify the existence of saddle point for the min-max problem  \eqref{mdp-5}. These conditions
define  a specific structure
of the ambiguity sets.
We begin with the first condition defined below.

 \begin{assumption}
\label{assum}
\begin{enumerate}
\item [{\rm (a)}]
 For every $s_t\in \s_t$ there is a kernel $P^*_{t}\in \PP_t$, potentially depending on $s_t$, with
\begin{equation}\label{saddle-b}
P^*_{t} \in
\argmax\limits_{P_t \in \PP_t} \tsum_{s_{t+1} \in \s_{t+1}} \tsum_{a_t \in \A(s_t)} P_t(s_{t+1} | s_t ,a_t) \mu_t(a_t) \big[
c_t(s_t, a_t ,s_{t+1})+V_{t+1}(s_{t+1}) \big],
\end{equation}
for   any $\mu_t \in \Delta_{|A_t(s_t)|}$.
\item [{\rm (b)}]
There is a kernel $P^*_{t}\in \PP_t$
  such that condition \eqref{saddle-b}  holds for all $s_t\in \s_t$  and any $\mu_t \in \Delta_{|A_t(s_t)|}$.
   \end{enumerate}
\end{assumption}

Since the set $\PP_t$ is assumed to be closed and hence  compact, a maximizer in the right hand side of \eqref{saddle-b} always exists. In condition (a) of the above assumption,  the kernel $P^*_{t}$ can be different for different values of the state $s_t\in \s_t$, while in condition (b) it is assumed to be the same (independent of $s_t$).
Clearly, condition (b) in the above assumption is stronger.

\begin{lemma}
\label{lem-cond}
Condition {\rm (b)} of Assumption {\rm \ref{assum}}  is equivalent to the following:
  there is a kernel $P^*_{t}\in \PP_t$ such that
\begin{equation}\label{saddle-a}
P^*_{t} \in
\argmax\limits_{P_t \in \PP_t} \tsum_{s_{t+1} \in \s_{t+1}} P_t(s_{t+1} | s_t ,a_t)   \big[
c_t(s_t, a_t ,s_{t+1})+V_{t+1}(s_{t+1}) \big], \;\forall
s_t\in \s_t,\;\forall a_t\in \A_t(s_t).
\end{equation}
Similar result holds for  condition {\rm (a)}.
\end{lemma}

\begin{proof}
If \eqref{saddle-b} holds for any probability $\mu_t$, then in particular it holds for Dirac measures and hence \eqref{saddle-a} follows. Conversely if \eqref{saddle-a}  holds for any $a_t\in \A_t(s_t)$, then it holds for convex combinations of such $a_t$, and hence \eqref{saddle-b} follows.
\end{proof}

Suppose that  Assumption \ref{assum}(a) holds   and consider
\begin{equation}\label{at}
\pi_t^*(\cdot|s_t)=\delta_{a_t^*(s_t)} \; \text{ with}\;a_t^*(s_t)\in\argmin_{a_t\in \A_t(s_t)} \tsum\limits_{s_{t+1} \in \s_{t+1}} P^*_{t}(s_{t+1} | s_t ,a_t)   \big[
c_t(s_t, a_t ,s_{t+1})+V_{t+1}(s_{t+1}) \big].
\end{equation}
   It follows directly from \eqref{saddle-b} and  \eqref{at}
that    $(\pi_t^* ,P_{t}^*)$ is a saddle point of the min-max problem  \eqref{mdp-5}.
Consequently  we have the following result.

\begin{theorem}
\label{th-duality}
Suppose that Assumption {\rm   \ref{assum}(a)} is fulfilled for all $t=1,...,T$,  and let $\pi_t^*(\cdot|s_t)$ be a minimizer defined in \eqref{at}. Then the following holds:
\begin{itemize}[itemsep=0.5pt, topsep=3pt]
\item[{\rm (i)}] $(\pi_t^*,P_{t}^*)$ is a saddle point of the min-max problem  \eqref{mdp-5}.
\item[{\rm (ii)}] There is no duality gap between  the primal problem \eqref{mdp-1} and its dual  \eqref{dual-1}.
\item[{\rm (iii)}]
$\{\pi^*_t\}$ is an optimal (non-randomized)  policy of the controller in the primal problem \eqref{mdp-1} and
$
 \gamma_t^*(s_t) = P_{t}^*$,  $ t = 1, \ldots, T,
 $
 is an optimal policy of the nature considered in the dual problem \eqref{dual-1}  (here $P_t^*$, viewed as an element of $\PP_t$,   can depend on $s_t\in \s_t$).
 \item[{\rm (iv)}]  If moreover  Assumption {\rm   \ref{assum}(b)} holds, then the above claims hold with $P_{t}^*$ independent of $s_t$.
 \end{itemize}
\end{theorem}

A sufficient condition ensuring Assumption \ref{assum}(b) is the so-called
$(\mathrm{s}, \mathrm{a})$-rectangularity.

\begin{definition}[$(\mathrm{s}, \mathrm{a})$-rectangularity,   \cite{iyen,Nilim2005}]
\label{def-rect}
Define
\begin{align}\label{sa_marginal}
\cP^t_{s_t,a_t}: = \{
q_t\in \Delta_{|\s_{t+1}|}: q_t(\cdot)  = P_t(\cdot | s_t,a_t)\; \text{for some}\;P_t \in \PP_t
\}.
\end{align}
The ambiguity set $\PP_t$ is said to be $(\mathrm{s}, \mathrm{a})$-rectangular if
\begin{equation}
\label{sarect}
\PP_t = \{
P_t: P_t(\cdot|s_t, a_t) \in \cP^t_{s_t,a_t}, \;\text{for all}\; s_t\in \s_t\;\text{and all}\;   a_t \in \A_t(s_t)\}.
\end{equation}
\end{definition}

If $\A_t$ is independent of $s_t$, then the
  $(\mathrm{s}, \mathrm{a})$-rectangularity  means  that $\PP_t$ can be represented as the direct product of $\cP^t_{s_t,a_t}$, $s_t\in \s_t$, $a_t\in \A_t$.

\begin{cor}\label{prop_sa_rec_strong_dual}
Suppose that the ambiguity sets $\PP_t$, $t=1,...,T$, are
$(\mathrm{s}, \mathrm{a})$-rectangular.
Then Assumption {\rm  \ref{assum}(b)} holds,
and conclusions {\rm (i) - (iv)} of Theorem {\rm \ref{th-duality}} follow.
\end{cor}

\begin{proof}
For $s_t\in \s_t$ and $a_t\in \A_t(s_t)$ consider the optimization problem
\begin{equation}\label{prob-1}
\max_{q_t\in \cP^t_{s_t,a_t}}\tsum_{s_{t+1}\in \s_{t+1}}q_t(s_{t+1})
 \big [c_t(s_t,a_t,s_{t+1})+V_{t+1}(s_{t+1})\big ].
\end{equation}
Since $\PP_t$ and hence  $\cP^t_{s_t,a_t}$ is compact, it follows that for any $s_t\in \s_t$ and $a_t\in \A_t(s_t)$ problem \eqref{prob-1} has an optimal solution $q^*_{t,s_t,a_t}$. Let $P_t^*(\cdot|s_t,a_t):=q^*_{t,s_t,a_t}(\cdot)$. By the $(\mathrm{s}, \mathrm{a})$-rectangularity of $\PP_t$ we have that $P^*_t\in \PP_t$, and
moreover   condition \eqref{saddle-a} holds for any $s_t\in \s_t$, and hence Assumption \ref{assum}(b) follows.
This completes the proof.
\end{proof}

\begin{remark}
\label{rem-dyn}
{\rm
Under the assumption of $(\mathrm{s},\mathrm{a})$-rectangularity,
dynamic equations of a   form similar to \eqref{mdp-5} have been established in \cite{iyen,Nilim2005}.
Nevertheless, it should be noted that the dynamic equations established there  are  for the static formulation to be discussed in Section \ref{subsec_static_form}. In comparison, Corollay \ref{prop_sa_rec_strong_dual} is established for the game formulation.
}  $\hfill \square$
\end{remark}

It turns out that   $(\mathrm{s}, \mathrm{a})$-rectangularity is not a necessary condition for Assumption \ref{assum}(b).
Indeed,  the following    simple two-stage problem ($T = 1$) gives an example  where the ambiguity set is not $(\mathrm{s}, \mathrm{a})$-rectangular while Assumption \ref{assum}(b) still holds.

\begin{figure}[t]
    \centering
    \includegraphics[width=0.38\textwidth]{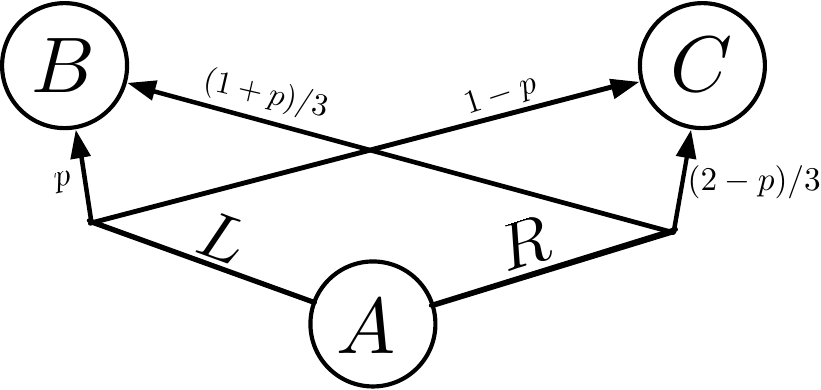} %
    \vspace{-0.1in}
    \caption{Two-stage distributionally robust MDP with a non-$(\mathrm{s}, \mathrm{a})$-rectangular ambiguity set. }
    \label{fig:example}
\end{figure}

\begin{example}\label{ex_example_non_necessary_of_sa}
{\rm
Consider a two-stage problem with three states $\{s_A, s_B, s_C\}$. Suppose $\s_1 = \{s_A \}$, $\s_2 = \{s_B, s_C\}$  and there  are two actions $a_L$ and $a_R$ at state $s_A$, i.e.   $\A_1=\{a_L,a_R\}$.
In addition, the cost functions are  given by $c_1 \equiv 0$, and $c_2(s_B) =v_B$ and $c_2(s_C)=v_C$ for some constants $v_B > v_C$.
In the considered setting,   for the corresponding  transition kernel $P_1$   we only need to specify the  probability  $P_1(s_B|s_A,a_L)$
 of moving to the state $s_B$ conditional on action $a_L$,   and
 the probability $P_1(s_B|s_A,a_R)$   conditional on action $a_R$,
since $P_1(s_C|s_A,a_L)=1-P_1(s_B|s_A,a_L)$ and $P_1(s_C|s_A,a_R)=1-P_1(s_B|s_A,a_R)$.
Let us consider the ambiguity set at stage $t=1$  defined as
\begin{equation}
\label{ambexam}
\PP_1 :=\big  \{
P_1: P_1(s_B|s_A, a_L) = p,  \;  P_1(s_B|s_A, a_R) =(1+p)/3,\;p\in [0,1]\big \}
\end{equation}
(see Figure \ref{fig:example}).
The corresponding set $\cP^1_{s_1,a_1}$, defined in \eqref{sa_marginal},  is determined  by a kernel  from the ambiguity set $\PP_1$.  It  follows   from \eqref{sarect} that the ambiguity set  $\PP_1$ is not $(\mathrm{s}, \mathrm{a})$-rectangular.

 On the other hand,   since $c_1\equiv 0$ and  $V_2(s_B)=c_2(s_B)$ is bigger than  $V_2(s_C)=c_2(s_C)$,   we have that  kernel  $P_1^*$   from the set $\PP_1$  for $p=1$, i.e.,
 defined as
\begin{align*}
P_1^*(s_B|s_A,a_L) = 1, ~ P_1^*(s_B |s_A, a_R) = \tfrac{2}{3},
\end{align*}
satisfies the condition of  Assumption \ref{assum}(b). $\hfill \square$
}
\end{example}

The above example falls into a broader class of ambiguity set defined below.

\begin{definition}[$\mathrm{r}$-rectangularity, \cite{goyal2023robust, goh2018data}]\label{def_r_rect}
The ambiguity set $\PP_t$ is said to be $\mathrm{r}$-rectangular if there exist
 sets $\W_1^t, \ldots, \W_m^t\subset \mathbb{R}^{\abs{\s_{t+1}}}$, and nonnegative  numbers  $\kappa^{s_t,a_t}_1,...,\kappa^{s_t,a_t}_m$ associated with each
$s_t \in \s_t$ and $a_t \in \A(s_t)$,
such that
\begin{align*}
\PP_t = \left\{
P_t: \exists w_i \in \W_i, i=1, \ldots, m, ~\text{and}~ P_t(\cdot|s_t, a_t) = \tsum_{i=1}^m \kappa^{s_t, a_t}_i w_i, ~\forall s_t \in \s_t, \forall a_t \in \A(s_t)
\right\}.
\end{align*}
\end{definition}

\begin{remark}\label{remark_r_rect}\
{\rm
The above definition of $\mathrm{r}$-rectangularity is slightly weaker than  the one considered in  \cite{goh2018data,goyal2023robust}.
Indeed, $\mathrm{r}$-rectangularity in  \cite{goyal2023robust, goh2018data}  requires additionally that  $\W_i \subseteq \Delta_{|\s_{t+1}|}$
and that  $\tsum_{i=1}^m \kappa^{s_t, a_t} _i= 1$. This additional condition is sufficient   for $\tsum_{i=1}^m  \kappa^{s_t, a_t}_i w_i$ to be a probability vector, although is not necessary.
Indeed,
the ambiguity set $\PP_1$ of   Example \ref{ex_example_non_necessary_of_sa} (defined in equation  \eqref{ambexam}) satisfies $\mathrm{r}$-rectangularity by taking $m=3$  and   $\kappa^{A, L}_1=1$, $\kappa^{A, L}_2=0$, $\kappa^{A, L}_3=1$, $\kappa^{A, R}_1=0$, $\kappa^{A, R}_2=1$, $\kappa^{A, R}_3= \tfrac{1}{3}$,  $\W_1=\{w_1\}$, $\W_2=\{w_2\}$,  $\W_3=\{w_3^p: p\in [0,1]\}$,  with $w_1=(0,1)\in \bbr^2$, $w_2=(\tfrac{1}{3}, \tfrac{2}{3})\in \bbr^2$ and $w_3^p=(p,-p)\in \bbr^2$.
 $\hfill \square$
}
\end{remark}

We proceed to verify Assumption \ref{assum}(b)  if the cost $c_t(s_t,a_t)$ does not depend on the next state $s_{t+1}$.\footnote{
Going forward whenever $\mathrm{r}$-rectangular set is involved, we assume the cost $c_t(s_t,a_t)$ does not involve the next state $s_{t+1}$.
}

\begin{cor}\label{prop_r_rectangular}
Suppose that the costs $c_t(s_t,a_t)$   do not depend on $s_{t+1}$, and the ambiguity sets $\PP_t$,  $t=1,...,T$,  are $\mathrm{r}$-rectangular.
Then Assumption {\rm  \ref{assum}(b)} holds,
and the conclusions {\rm (i)} -  {\rm (iv)}
of Theorem {\rm \ref{th-duality}} follow.
\end{cor}

\begin{proof}
For $\mathrm{r}$-rectangular sets the maximization problem in \eqref{saddle-b} is equivalent to
\begin{eqnarray*}
 &\max\limits_{P_t \in \PP_t}&\!\!\! \tsum_{s_{t+1} \in \s_{t+1}} \tsum_{a_t \in \A(s_t)} P_t(s_{t+1} | s_t ,a_t) \mu_t(a_t) \big[
c_t(s_t, a_t )+V_{t+1}(s_{t+1}) \big]  \\
  &&=  \max_{w_i \in \W_i^t, 1\leq i \leq m}  \left\{
\tsum_{a_t \in \A(s_t)} \mu_t(a_t) c_t(s_t, a_t)
+ \tsum_{a_t \in \A(s_t)} \mu_t(a_t) \Big(
\tsum_{i=1}^m \kappa^{s_t, a_t}_i \cdot w_i^\top V_{t+1}
\Big)
\right\},
\end{eqnarray*}
where we denote by  $V_{t+1}$  the $|\s_{ t+1}|$-dimensional vector with components $V_{t+1}(s_{ t+1})$, $s_{ t+1}\in \s_{ t+1}$,  and $a^\top b=\tsum_{i=1}^n a_i b_i$ is the scalar product of vectors $a,b\in \bbr^n$.

Since $\kappa^{s_t, a_t}_i \geq 0$,
  in order to satisfy Assumption \ref{assum}(b)  it suffices to choose $P_t^*$ defined as
\begin{equation}
P_t^*(\cdot |s_t, a_t) :=  \tsum_{i=1}^m \kappa^{s_t, a_t}_i w_i^*,
\end{equation}
where $w_i^* \in \argmax\limits_{w_i \in \W^t_i} w_i^\top V_{t+1}$.
\end{proof}

\vspace{0.1in}
 In the absence of the first condition (Assumption \ref{assum}), we introduce the following alternative condition that certifies the strong duality of the game formulation.

\begin{assumption}\label{assump_margin_convex}
For all $t=1,...,T$,  and $s_t\in \s_t$,
the state-wise marginalization of $\PP_t$, defined as
\begin{align}\label{eq_marginal_s}
\cP^t_{s_t}: = \{
q_t\in \Delta_{|\s_{t+1}\times \A_t(s_t)|}: q_t (\cdot,\cdot)= P_t(\cdot|s_t, \cdot) ~ \text{for some} ~ P_t \in \PP_t
\},
\end{align}
is convex.
\end{assumption}

 It can be noted that the above assumption holds if the set $\PP_t$ is convex.

\begin{theorem}\label{prop_s_rect_convex_marginal}
Suppose Assumption {\rm  \ref{assump_margin_convex}} is fulfilled.
Then the following holds.
 \begin{itemize}[itemsep=0.5pt, topsep=3pt]
\item[{\rm (i)}] There exists a saddle point $(\pi_t^*,P_{t}^*)$ for the min-max problem  \eqref{mdp-5}, and hence
there is no duality gap between  the primal problem \eqref{mdp-1} and its dual  \eqref{dual-1}.
\item[{\rm (ii)}]  The primal problem \eqref{mdp-1} has a non-randomized optimal policy iff the min-max problem
\eqref{mark-2}  has a saddle point for all $t=1,...,T$,  and $s_t\in \s_t$.
 \end{itemize}
\end{theorem}

\begin{proof}
Clearly, \eqref{mdp-5} is equivalent to the following optimization problem
\begin{equation}
\label{probminm}
\min_{\pi_t(\cdot|s_t) \in \Delta_{|\A_t(s_t)|}} \sup_{P_t(\cdot|s_t, \cdot) \in \cP_{s_t}^t } \tsum\limits_{s_{t+1} \in \s_{t+1}} \tsum\limits_{a_t \in \A_t(s_t)} P_t(s_{t+1} |  s_t, a_t) \pi_t(a_t | s_t) \big[
c_t(s_t, a_t ,s_{t+1})+V_{t+1}(s_{t+1}) \big].
\end{equation}
The objective function of problem \eqref{probminm} is linear with respect to $\pi_t(\cdot|s_t)$ and linear with respect to $P_t(\cdot|s_t, \cdot)$, and both sets $\Delta_{|\A_t(s_t)|}$ and $\cP_{s_t}^t $ are convex and the set $\Delta_{|\A_t(s_t)|}$ is compact. Therefore the strong duality follows by Sion's theorem.
  Moreover the primal and dual problems posses optimal solutions and hence there exists the corresponding saddle point.
  This completes the proof of (i).

 Now let us observe that the  minimum in the max-min problem \eqref{dual-4} is always attained at Dirac measure, and hence
  we can write
   \begin{equation}\label{dual-4non}
Q_t(s_t)=\max_{P_t\in \PP_t} \min_{a_t \in  \A_t(s_t)}   \tsum\limits_{s_{t+1} \in \s_{t+1}}  P_t(s_{t+1} | s_t ,a_t)  \big[
c_t(s_t, a_t ,s_{t+1})+Q_{t+1}(s_{t+1}) \big].
\end{equation}
 Suppose that \eqref{mark-2}  has a saddle point $(a_t^*,P_t^*)\in \A_t(s_t)\times \PP_t$. Then $(a_t^*,P_t^*)$ is a saddle point of
 \eqref{dual-4non}, and hence $\{\delta_{a^*_t(s_t)}\}$ is
a  non-randomized optimal policy of the controller.
Conversely suppose that the primal problem has a non-randomized policy. By (i) we have   that there is no duality gap between the primal and dual problems. Therefore  there is  saddle point $(\pi_t^*,P_{t}^*)$ with $\pi(\cdot|s_t)$ being Dirac measure.
This completes the proof of (ii).
\end{proof}

We next introduce an example demonstrating the essential role of convexity in Assumption \ref{assump_margin_convex}.

\begin{example}\label{ex_s_rect_non_saddle}
{\rm
Consider the same MDP as defined in Example \ref{ex_example_non_necessary_of_sa}, except that we have
\begin{align*}
\PP_1 =\big  \{
P_1: P_1(s_B|s_A,a_L) = p, ~ P_1(s_B|s_A, a_R) = 1- p, \;  p \in [0,\tfrac{1}{4}] \cup [\tfrac{3}{4}, 1]  \}.
\end{align*}
Clearly $\cP^1_{s_A}$ is not convex.
Direct computation shows
\begin{align}
V_1(s_A) & = \min_{\mu \in \Delta_{\{L, R\}}} \max_{P_1 \in \PP_1}
\tsum_{s_{2} \in \s_{2}} \tsum_{a_1 \in \{a_L, a_R\}} P_1(s_2 |s_A ,a_1) \mu(a_1)
c_2(s_2)
= \tfrac{v_B + v_C}{2} ,  \nonumber \\
Q_1(s_A) & = \max_{P_1 \in \PP_1}  \min_{\mu \in \Delta_{\{L, R\}}}
\tsum_{s_{2} \in \s_{2}} \tsum_{a_1 \in \{a_L, a_R\}} P_1(s_2 | s_A ,a_1) \mu(a_1)
c_2(s_2)
=  \tfrac{1}{4} v_B + \tfrac{3}{4} v_C, \label{eq_strick_weak_dual}
\end{align}
and hence $Q_1(s_A)< V_1(s_A)$.
$\hfill \square$
}
\end{example}

As shown by Example \ref{ex_s_rect_non_saddle},  without the convexity in Assumption \ref{assump_margin_convex}, the weak duality \eqref{opt-1} can be  strict.
At this point, it might be worthing noting  that with Assumption \ref{assump_margin_convex} the optimal policies of the controller  can be all
randomized, i.e.,  if the necessary condition specified in Theorem \ref{prop_s_rect_convex_marginal}(ii) does not hold.
The next example illustrates the aforementioned observation. 
A similar claim is made in \cite{kuhn2013} for the static formulation we discuss in Section \ref{subsec_static_form}.

\begin{example}\label{example_randomized_opt_policy}
{\rm
Consider the same MDP as defined in Example \ref{ex_s_rect_non_saddle}, except that in the definition of $\PP_1$ we let $p \in [0,1]$.
Direct computation shows that when $\mu = \delta_{a_{_L}}$, the unique solution of   \eqref{saddle-b} is given by
$P_1(s_B | s_A, a_L) = 1$ and $P_1(s_B| s_A, a_R) = 0$.
In comparison, when $\mu = \delta_{a_R}$, the unique solution of \eqref{saddle-b} is given by
$P_1(s_B | s_A, a_L) = 0$ and $P_1(s_B| s_A, a_R) = 1$.
Consequently Assumption \ref{assum}(a) does not hold.
Nevertheless,  $\PP^1_{s_A}$ is convex and closed, and hence Theorem \ref{prop_s_rect_convex_marginal} applies and strong duality holds.
Finally, it can be readily verified that the optimal policy $\pi^*$ of \eqref{mdp-1} is unique and is given by
\begin{align*}
\pi^*_1(a_L|s_A) = \pi^*_1(a_R|s_A) = 1/2.
\end{align*}
Note that the policy  $\pi^*$ is randomized.
 $\hfill \square$
}
\end{example}

Example \ref{example_randomized_opt_policy} should be contrasted with conclusion (iii) of Theorem \ref{th-duality}, which always certifies the existence of non-randomized optimal policies in the presence of Assumption \ref{assum}.

\begin{remark}[Nature's Randomized Policies]
{\rm
It is also possible to consider randomized policies of the nature when the kernel $P_t$ is chosen at random according to a probability distribution supported on the ambiguity set $\PP_t$.
Then the counterpart of dynamic equations \eqref{mdp-5} can be written as
\begin{equation}\label{counterpart}
V_t(s_t)=\min_{\pi_t(\cdot|s_t) \in \Delta_{|\A_t(s_t)|}} \sup_{D_t \in \cD_t} \bbe_{P_t \sim D_t} \bbe^{\pi_t, P_t}[
c_t(s_t, a_t ,s_{t+1})+V_{t+1}(s_{t+1})],
\end{equation}
where $\cD_t$ is the set of all probability distributions on $\PP_t$ equipped with its Borel sigma algebra.  By Sion's theorem the min-max  problem \eqref{counterpart} has a saddle point (note that here the existence of saddle point does not require the convexity of $\cP^t_{s_t}$). Using an argument similar to the proof of  part (ii) of  Theorem \ref{prop_s_rect_convex_marginal}, we have here that the nature possesses a non-randomized policy iff the min-max problem \eqref{mdp-5} has a saddle point for all $t=1,...,T,$ and $s_t\in \s_t$.
} $\hfill \square$
\end{remark}

\subsection{History-dependent policies}
\label{histdepen}

In this section we discuss history-dependent policies of the controller and/or  the nature. We denote by $\widehat{\Pi}$ and
$\widehat{\Gamma}$ the sets of history-dependent policies of the controller and    the nature, respectively, and consider the counterparts of the primal problem  \eqref{mdp-1}  and its dual \eqref{dual-1} when $\Pi$ and $\Gamma$ are replaced by their history-dependent counterparts $\widehat{\Pi}$ and
$\widehat{\Gamma}$. Recall that $\ch_t$ denotes history of the decision process (see~\eqref{history}).

 \begin{proposition}\label{prop_hist}
 The following holds. {\rm (i)} Suppose the controller is Markovian and the nature is history-dependent, i.e., consider the problem
 \begin{equation}\label{mdp-hist}
 \min_{\pi\in \Pi}\max_{\gamma\in \widehat{\Gamma}}\bbe^{\pi,\gamma}
 \left [ \tsum_{t=1}^{T}
c_t(s_t,a_t,s_{t+1})+c_{T+1}(s_{T+1})
\right].
\end{equation}
 Then the dynamic equations \eqref{mdp-3} and
 \eqref{mdp-5} still hold.

 {\rm (ii)} Suppose both the controller and the nature are history-dependent, i.e.,  consider the problem
\begin{equation}\label{mdp-his2}
 \min_{\pi\in \widehat{\Pi}}\max_{\gamma\in \widehat{\Gamma}}\bbe^{\pi,\gamma}
 \left [ \tsum_{t=1}^{T}
c_t(s_t,a_t,s_{t+1})+c_{T+1}(s_{T+1})
\right].
\end{equation}
Then for any $\pi\in \widehat{\Pi}$, we have
$\V^\pi_{\ch_T} (s_{T+1}) = c(s_{T+1})$, and for $t=T,...,1$,
 \begin{align}\label{mdp-3_hist}
\V^\pi_{\ch_{t-1}}(s_t)  = \max_{P_t \in \PP_t}  {\tsum_{s_{t+1} \in \s_{t+1}}  \tsum_{a_t \in \A_t(s_t)} P_t(s_{t+1} | s_t ,a_t) \pi_t(a_t | \ch_{t-1}, s_t) \big[
c_t(s_t, a_t ,s_{t+1})+\V^\pi_{\ch_{t}}(s_{t+1}) \big]},
\end{align}
where $\V^\pi_{\ch_{t-1}}(s_t) $ is equal  to  the optimal value of problem
\begin{align}\label{eq_inner_max_starting_t_hist}
\max_{\gamma \in \widehat{\Gamma}}  \bbe^{\pi, \gamma}\left [ \tsum_{i=t}^{T}
c_i(s_i,a_i,s_{i+1})+c_{T+1}(s_{T+1})
\large| \ch_{t-1}, s_t
\right].
\end{align}
    Moreover,
 $V_t(s_t)$ defined in \eqref{mdp-5},  is equal  to the optimal value of
 \begin{align}\label{eq_primal_starting_at_t_hist}
 \min_{\pi \in \widehat{\Pi}}
 \max_{\gamma \in \widehat{\Gamma}}  \bbe^{\pi, \gamma}\left [ \tsum_{i=t}^{T}
c_i(s_i,a_i,s_{i+1})+c_{T+1}(s_{T+1}) \big| \ch_{t-1} ,s_t
\right].
 \end{align}
\end{proposition}

\begin{proof}
For claim {\rm (i)}, note that the argument of Proposition \ref{pr-dyncontr} holds even if the nature uses history-dependent policies, since  optimizing over Markovian policies is the same as optimizing over history-dependent policies for regular MDPs \cite{puterman2014markov}.
Consequently Proposition \ref{prop_game_primal_dp} also holds.

We proceed to establish  claim {\rm (ii)}.
Let us first consider \eqref{eq_inner_max_starting_t_hist}.
Fixing any $\gamma \in \widehat{\Gamma}$, it is immediate that $\V^{\pi,\gamma}_{\ch_T} (s_{T+1}) = c(s_{T+1})$,
\begin{align}\label{eq_dp_fix_nature_policy_hist}
\V^{\pi, \gamma}_{\ch_{t-1}}(s_t)  =  {\tsum_{s_{t+1} \in \s_{t+1}}  \tsum_{a_t \in \A_t(s_t)} P_t(s_{t+1} | s_t ,a_t) \pi_t(a_t | \ch_{t-1}, s_t) \big[
c_t(s_t, a_t ,s_{t+1})+\V^{\pi, \gamma}_{\ch_{t}}(s_{t+1}) \big]}
\end{align}
defines $\{\V^{\pi, \gamma}_{\ch_{t-1}}\}$ such that
\begin{align*}
\V^{\pi, \gamma}_{\ch_{t-1}}(s_t) = \bbe^{\pi, \gamma}\left [ \tsum_{i=t}^{T}
c_i(s_i,a_i,s_{i+1})+c_{T+1}(s_{T+1})
\large| \ch_{t-1}, s_t
\right].
\end{align*}
Here $P_t$ in \eqref{eq_dp_fix_nature_policy_hist} is specified by nature's policy $\gamma$, i.e., $P_t = \gamma_t (\ch_{t-1}, s_t)$.
We now show by induction that $V^\pi_{\ch_{t}} \geq \V^{\pi, \gamma}_{\ch_{t}}$ for $t = T, \ldots, 0$.
The claim is obvious at $t = T$. In addition,
\begin{align*}
\V^\pi_{\ch_{t-1}}(s_t)  & \overset{(a)}{\geq} \max_{P_t' \in \PP_t}  {\tsum_{s_{t+1} \in \s_{t+1}}  \tsum_{a_t \in \A_t(s_t)} P_t'(s_{t+1} | s_t ,a_t) \pi_t(a_t | \ch_{t-1}, s_t) \big[
c_t(s_t, a_t ,s_{t+1})+\V^{\pi,\gamma}_{\ch_{t}}(s_{t+1}) \big]} \\
& {\geq} {\tsum_{s_{t+1} \in \s_{t+1}}  \tsum_{a_t \in \A_t(s_t)} P_t(s_{t+1} | s_t ,a_t) \pi_t(a_t | \ch_{t-1}, s_t) \big[
c_t(s_t, a_t ,s_{t+1})+\V^{\pi,\gamma}_{\ch_{t}}(s_{t+1}) \big]} \\
& \overset{(b)}{=} \V^{\pi, \gamma}_{\ch_{t-1}}(s_t),
\end{align*}
from which we complete the induction step.
Here $(a)$ follows from the induction hypothesis,
and $(b)$ follows from \eqref{eq_dp_fix_nature_policy_hist}.

Now let $\tilde{\gamma}$ be the history-dependent policy of the nature that attains the maximum in \eqref{mdp-3_hist},
then from \eqref{eq_dp_fix_nature_policy_hist} it holds that
$V^\pi_{\ch_{t}} = \V^{\pi, \tilde{\gamma}}_{\ch_{t}}$.
Consequently, we obtain $V^\pi_{\ch_{t}}  = \max_{\gamma}  \V^{\pi, \gamma}_{\ch_{t}}$ and \eqref{eq_inner_max_starting_t_hist}.

With \eqref{mdp-3_hist} in place, one can follow essentially the same lines as in the proof of Proposition \ref{prop_game_primal_dp} and readily obtain that dynamic equations  $V_{\ch_T}(s_{T+1})=  c_{T+1}(s_{T+1})$,
\begin{equation}\label{mdp-5_hist}
V_{\ch_{t-1}}(s_t)=\min_{\pi_t(\cdot|\ch_{t-1}, s_t) \in \Delta_{|\A_t(s_t)|}} \max_{P_t \in \PP_t} \tsum\limits_{s_{t+1} \in \s_{t+1}} \tsum\limits_{a_t \in \A_t(s_t)} P_t(s_{t+1} | s_t ,a_t) \pi_t(a_t | \ch_{t-1} , s_t) \big[
c_t(s_t, a_t ,s_{t+1})+V_{\ch_t}(s_{t+1}) \big],
\end{equation}
define  $V_{\ch_{t-1}}(s_t)$,  which equals the optimal value of \eqref{eq_primal_starting_at_t_hist}.
 To conclude our proof for claim {\rm (ii)}, it suffices to show that $V_{\ch_{t}}(\cdot)$ only depends on $t+1$ and does not depend on $\ch_{t}$.
 We now show this by induction.
Clearly, this holds trivially at $t = T$.
Suppose the claim holds at stage $t$, then
from \eqref{mdp-5_hist} we obtain
\begin{equation*}
V_{\ch_{t-1}}(s_t)=\min_{q_t  \in \Delta_{|\A_t(s_t)|}} \max_{P_t \in \PP_t} \tsum\limits_{s_{t+1} \in \s_{t+1}} \tsum\limits_{a_t \in \A_t(s_t)} P_t(s_{t+1} | s_t ,a_t) q_t(a_t) \big[
c_t(s_t, a_t ,s_{t+1})+V_{t+1}(s_{t+1}) \big],
\end{equation*}
from which it is clear that $V_{\ch_{t-1}}(\cdot)$ only depends on $t$ but not $\ch_{t-1}$.
From this we complete the induction step and the proof for claim {\rm (ii)} is completed.
\end{proof}

In  \cite[section 5]{WangBlanchet2023}  examples are presented when the controller is allowed  history-dependent policies while the nature is  Markovian,  which are not amendable to writing   dynamic programming equations similar to \eqref{mdp-5}.
We show now that  the respective dynamic equations still hold provided there is no duality gap between the primal problem \eqref{mdp-1} and its dual \eqref{dual-1}.

\begin{proposition}
Consider counterpart of \eqref{mdp-1} when the controller uses history-dependent policies $\widehat{\Pi}$ and the nature uses Markovian policies $\Gamma$, i.e.,
\begin{align}\label{prob_hist_controller_markov_nature}
 \min_{\pi\in \widehat{\Pi}}\max_{\gamma\in {\Gamma}} \bbe^{\pi,\gamma}
 \left [ \tsum_{t=1}^{T}
c_t(s_t,a_t,s_{t+1})+c_{T+1}(s_{T+1})
\right].
\end{align}
If $ \mathrm{OPT}(\ref{mdp-1}) =  \mathrm{OPT}(\ref{dual-1}) $ holds, then there exists an optimal policy of \eqref{prob_hist_controller_markov_nature} that is Markovian,
and dynamic equations \eqref{mdp-5} still hold for problem \eqref{prob_hist_controller_markov_nature}.
\end{proposition}

\begin{proof}
Let us denote $f(\pi, \gamma) := \bbe^{\pi,\gamma}
  [ \tsum_{t=1}^{T}
c_t(s_t,a_t,s_{t+1})+c_{T+1}(s_{T+1})
]$.
Then the following is immediate:
\begin{align*}
\max_{\gamma \in \Gamma} \min_{\pi \in \Pi} f(\pi, \gamma) \overset{(a)}{=}
\max_{\gamma \in \Gamma} \min_{\pi \in \widehat{\Pi}} f(\pi, \gamma)
\leq
\min_{\pi \in \widehat{\Pi}}  \max_{\gamma \in \Gamma} f(\pi, \gamma)
\leq
\min_{\pi \in {\Pi}}  \max_{\gamma \in \Gamma} f(\pi, \gamma) ,
\end{align*}
where $(a)$ follows from the fact that $\min_{\pi \in \widehat{\Pi}} f(\pi, \gamma)$ is a regular MDP of the controller as the nature is Markovian,
subsequently for the controller minimizing over Markovian policies is equivalent to minimizing over history-dependent policies.
If $  \mathrm{OPT}(\ref{mdp-1}) = \mathrm{OPT}(\ref{dual-1}) $ holds,
then clearly all the above inequalities indeed hold with equality, from which we conclude the proof.
\end{proof}

Together with Theorem \ref{prop_s_rect_convex_marginal}(i)  this implies the following.

\begin{cor}\label{cor_existence_markov_opt_policy}
Suppose Assumption {\rm   \ref{assum}(a)} or Assumption {\rm  \ref{assump_margin_convex}} is fulfilled. Then there exists an optimal policy of the controller in \eqref{prob_hist_controller_markov_nature} that is Markovian.
In particular, this holds when the ambiguity set $\PP_t$ is $\mathrm{(s,a)}$- or $\mathrm{r}$-rectangular, or has convex state-wise marginalization (defined in \eqref{eq_marginal_s}).
\end{cor}

In Section \ref{subsec_static_form}, we will introduce $\mathrm{sr}$-rectangular sets.
One can readily check that conclusion of Corollary \ref{cor_existence_markov_opt_policy} also holds for any $\mathrm{sr}$-rectangular set if the set is $\mathrm{s}$-convex (Definition \ref{def_sr_rectangular}).

\setcounter{equation}{0}
\section{Static formulation of distributionally robust MDP}\label{subsec_static_form}
Let us consider the following  static formulation of distributionally robust MDP,
\begin{align}\label{eq_static_formulation}
 \min_{\pi\in \Pi}\sup_{P_1\in \PP_1, \ldots, P_T \in \PP_T} \bbe^{\pi,\{P_t\}_{t = 1}^T}
 \left [ \tsum_{t=1}^{T}
c_t(s_t,a_t,s_{t+1})+c_{T+1}(s_{T+1})
\right].
\end{align}
As it was already mentioned,    in some publications   the static formulation     is  referred to as {\em robust } MDPs.
The essential difference between \eqref{eq_static_formulation} and the game formulation \eqref{mdp-1} is that in \eqref{eq_static_formulation} the nature  chooses the kernels  prior to the realization of the Markov process.
The   dual  of \eqref{eq_static_formulation} is
\begin{align}\label{eq_static_formulation_dual}
\sup_{P_1\in \PP_1, \ldots, P_T \in \PP_T} \min_{\pi\in \Pi} \bbe^{\pi,\{P_t\}_{t = 1}^T}
 \left [ \tsum_{t=1}^{T}
c_t(s_t,a_t,s_{t+1})+c_{T+1}(s_{T+1})
\right].
\end{align}
Formulations \eqref{eq_static_formulation} and \eqref{eq_static_formulation_dual} can be viewed as
 the static counterparts of the respective problems  \eqref{mdp-1} and \eqref{dual-1}.

Below we first discuss when the static formulation (for both primal and dual problems) is equivalent to the game formulation discussed in Section \ref{subsec_game_form}.

\begin{theorem}\label{prop_form_equivalence}
The following holds.
{\rm (i)}  If for any policy $\pi\in \Pi$,  there exists a solution $P_{t}^*$ of \eqref{mdp-3} that is independent of
 $s_t$, for every $s_t \in \s_t$ and $t = 1,\ldots, T$,
then the primal problem \eqref{mdp-1} of the game  and its static counterpart    \eqref{eq_static_formulation}   are equivalent, i.e.,
\begin{align}
 \mathrm{OPT} \eqref{mdp-1} =  \mathrm{OPT} \eqref{eq_static_formulation},
  \end{align}
  and the dynamic equations \eqref{mdp-5} hold for the static problem  \eqref{eq_static_formulation} as well.

{\rm (ii)}  If for problem  \eqref{dual-4} (not necessarily having a saddle point), there exits a solution $P_{t}^*$ such that $P_{t}^*$ is independent of
 $s_t$, for every $s_t \in \s_t$ and $t = 1,\ldots, T$,
then the dual problem   \eqref{dual-1} of the game and its  static counterpart   \eqref{eq_static_formulation_dual}    are equivalent, i.e.,
\begin{align}
   \mathrm{OPT} \eqref{dual-1} = \mathrm{OPT} \eqref{eq_static_formulation_dual} ,
\end{align}
and the dynamic equations \eqref{dual-4} hold for the static problem  \eqref{eq_static_formulation_dual} as well.
\end{theorem}

\begin{proof}
For the   claim (i), given the specified   condition, it can be readily seen that for any $\pi\in \Pi$,
\begin{equation}
\label{eq-ineq}
\begin{array}{lll}
    \max\limits_{\gamma \in \Gamma}&  \bbe^{\pi, \gamma} \left [ \tsum_{t=1}^{T}
c_t(s_t,a_t,s_{t+1})+c_{T+1}(s_{T+1})
\right]
 \\
  &=     \bbe^{\pi, \{P_t^*\}_{t=1}^T}  \left [ \tsum_{t=1}^{T}
c_t(s_t,a_t,s_{t+1})+c_{T+1}(s_{T+1})
\right]
\\
  & \leq   \sup\limits_{P_1\in \PP_1, \ldots, P_T \in \PP_T} \bbe^{\pi,\{P_t\}_{t = 1}^T}
 \left [ \tsum_{t=1}^{T}
c_t(s_t,a_t,s_{t+1})+c_{T+1}(s_{T+1})
\right].
\end{array}
\end{equation}
On the other hand, it holds trivially that the first term in \eqref{eq-ineq} is at least that of the last term in the above relation, which implies that the above inequality holds with equality. Since this holds for any $\pi\in \Pi$, the   claim (i) follows.

For the  claim (ii), from the specified  condition we have that for the dual problem \eqref{dual-1}, the optimal policy is given by $\gamma_t^*(s_t) = P_t^*$, and consequently
\begin{align*}
\mathrm{OPT} \eqref{dual-1}  & = \min_{\pi\in \Pi} \bbe^{\pi,\gamma^*}
 \left [ \tsum_{t=1}^{T}
r_t(s_t,a_t,s_{t+1})+r_{T+1}(s_{T+1})
\right]  \\
&  = \min_{\pi\in \Pi} \bbe^{\pi,\{P_t^*\}_{t = 1}^T}
 \left [ \tsum_{t=1}^{T}
c_t(s_t,a_t,s_{t+1})+c_{T+1}(s_{T+1})
\right]
\leq
\mathrm{OPT} \eqref{eq_static_formulation_dual} .
\end{align*}
Combining the above relation with $\mathrm{OPT} \eqref{eq_static_formulation_dual} \leq \mathrm{OPT} \eqref{dual-1}$ yields the claim (ii).
\end{proof}

Clearly, the condition in Theorem \ref{prop_form_equivalence} hold for $(\mathrm{s},\mathrm{a})$- and $\mathrm{r}$-rectangular sets discussed in Section \ref{subsec_game_form}.
Below we introduce another commonly used ambiguity set satisfying this condition (the verification uses a similar argument as in Corollary \ref{prop_sa_rec_strong_dual}).

\begin{definition}[$\mathrm{s}$-rectangularity, \cite{Tallec, kuhn2013}]\label{def_s_rect}
The ambiguity set $\PP_t$ is said to be $\mathrm{s}$-rectangular if the following holds:
\begin{align*}
\PP_t = \{
P_t: P_t(\cdot|s_t, \cdot) \in \cP^t_{s_t}, ~ \forall s_t\in \s_t
\},
\end{align*}
where $\cP^t_{s_t}$ is defined as in \eqref{eq_marginal_s}.
\end{definition}

 Theorem \ref{prop_form_equivalence} allows one to immediately obtain the dynamic equations of the static formulation for the ambiguity sets commonly used for the static formulation.

\begin{cor}\label{cor_form_equivalence}
For  the $(\mathrm{s},\mathrm{a})$-, $\mathrm{s}$- and $\mathrm{r}$-rectangular sets, the game formulation is equivalent to the static formulation. In particular, \eqref{mdp-5} provides the dynamic equations of the primal problem \eqref{eq_static_formulation} of the static formulation,
and \eqref{dual-4} provides the dynamic equations of the dual problem \eqref{eq_static_formulation_dual} of the static formulation.
\end{cor}

\begin{remark}\label{remark_form_equivalence}
{\rm
It might be worth mentioning that prior studies of the static formulation adopt case-by-case analyses in establishing the corresponding dynamic equations for different rectangular sets (cf., \cite{kuhn2013, goyal2023robust, iyen,Nilim2005}).
The formulation equivalence we establish here allows one to provide a unified treatment for different ambiguity sets with a rather simple approach.
}
$\hfill \square$
\end{remark}

 Theorem \ref{prop_form_equivalence} also provides a constructive way to reformulate a game formulation with potentially non-rectangular ambiguity set into the static formulation with a rectangular one. Recall that $\cP^t_{s_t}$, defined in  \eqref{eq_marginal_s},  denotes the state-wise marginalization of $\PP_t$.

\begin{theorem}\label{cor_reform_game}
The game formulation \eqref{mdp-1} (resp. \eqref{dual-1}) with an arbitrary ambiguity set $\PP_t$ is equivalent to another game with
a (possibly enlarged) $\mathrm{s}$-rectangular ambiguity set $\PP'_t$, where
\begin{equation}\label{largeambig}
\PP_t' =  \{
P_t: P_t(\cdot|s_t, \cdot) \in \cP^t_{s_t},  \forall s_t \in \s_t
\}
\end{equation}
Consequently, the game formulation \eqref{mdp-1} (resp. \eqref{dual-1})   can   be equivalently reformulated into
the static formulation \eqref{eq_static_formulation} (resp. \eqref{eq_static_formulation_dual}) with $\mathrm{s}$-rectangular ambiguity set $\PP'_t$.
\end{theorem}

\begin{proof}
Consider the dynamic equations \eqref{mdp-5} (resp. \eqref{dual-4})  for the primal problem \eqref{mdp-1} (resp. for the dual problem  \eqref{dual-1}). Replacing the ambiguity set $\PP_t$ with the set $\PP'_t$ does not change the respective maximum.
The claim   follows by   invoking Corollary \ref{cor_form_equivalence}.
\end{proof}

\begin{remark}
{\rm
For the static formulation \eqref{eq_static_formulation}, identifying ambiguity sets beyond $\mathrm{s}$-rectangularity while maintaining favorable dynamic equations has been an active line of research.
For instance, it is shown in \cite{goyal2023robust} that the class of $\mathrm{r}$-rectangular sets has elements that can not be represented by any $\mathrm{s}$-rectangular set.
On the other hand, in view of Theorem \ref{cor_reform_game}, for formulation purposes there seems to be limited motivation in seeking broader classes of rectangular sets.
 As an illustration, we have the following simple observation showing that any $\mathrm{r}$-rectangular distributionally robust MDP is equivalent to a properly constructed $\mathrm{s}$-rectangular counterpart.\footnote{
For practical purposes, a more delicate construction of ambiguity set can be still favorable when data or domain knowledge can be used within such a construction. The enlarged ambiguity set could also be significantly larger than the original one, leading to difficulties from a computational perspective.
We make an attempt at this direction and construct in Section \ref{subsec_sr_rectangular} a richer class of ambiguity sets that encompass $\mathrm{(s,a)}$-, $\mathrm{s}$-, and $\mathrm{r}$-rectangular sets.
}
}
$\hfill \square$
 \end{remark}

\begin{cor}\label{cor_from_r_to_s}
Any static formulation \eqref{eq_static_formulation} (resp. \eqref{eq_static_formulation_dual}) with an $\mathrm{r}$-rectangular ambiguity set $\PP_t$ is equivalent to
another static formulation with an (enlarged)  $\mathrm{s}$-rectangular ambiguity set $\PP_t'$.
\end{cor}
\begin{proof}
The claim follows by noting that \eqref{eq_static_formulation} with an $\mathrm{r}$-rectangular ambiguity set is equivalent to the game formulation \eqref{mdp-1} and subsequently invoking Theorem \ref{cor_reform_game}.
\end{proof}

It appears that $\mathrm{r}$-rectangular sets are believed to be favorable over $\mathrm{s}$-rectangular sets, as  conceptually it allows coupling of transition probabilities across different states, and thus can produce less conservative policy than  $\mathrm{s}$-rectangular sets \cite{goh2018data, goyal2023robust}.
Corollary \ref{cor_from_r_to_s} above shows that one can obtain the same robust policy associated with any  $\mathrm{r}$-rectangular set by using a properly constructed $\mathrm{s}$-rectangualr set.

We are ready to discuss the duality of the static formulation.

\begin{theorem}\label{prop_form_eq_and_strong_static_dual}
Suppose for every $t = 1,\ldots T$ and $s_t \in \s_t$, a saddle point $(\pi_t^*, P_{t}^*)$ exists for \eqref{mdp-5}, and $P_{t}^*$ does not depend on $s_t$, then
\begin{align}\label{equivalence_and_strong_dual}
 \mathrm{OPT} \eqref{mdp-1} =  \mathrm{OPT} \eqref{dual-1}  =  \mathrm{OPT} \eqref{eq_static_formulation} = \mathrm{OPT} \eqref{eq_static_formulation_dual}.
\end{align}
That is,  the static formulation (for both primal and dual problems) is equivalent to the game formulation,
and the strong duality holds for the static formulation.
\end{theorem}

\begin{proof}
From Remark \ref{remark_strong_dual_general}, we have $\mathrm{OPT} \eqref{mdp-1} =  \mathrm{OPT} \eqref{dual-1}$.
Consequently a saddle point of \eqref{mdp-5} is also a saddle point of \eqref{dual-4}.
Since $P_{t}^*$ does not depend on $s_t$, we can invoke the second claim of Theorem \ref{prop_form_equivalence} and obtain that $ \mathrm{OPT} \eqref{dual-1} = \mathrm{OPT} \eqref{eq_static_formulation_dual}$.
The desired claim then follows from the additional trivial observations:
$
 \mathrm{OPT} \eqref{eq_static_formulation_dual} \leq  \mathrm{OPT} \eqref{eq_static_formulation} \leq \mathrm{OPT} \eqref{mdp-1}
$
\end{proof}

\begin{cor}
\label{cor-equality}
 {\rm (i)}
Under Assumption {\rm \ref{assum}(b)}, the conclusions of Theorem {\rm \ref{prop_form_eq_and_strong_static_dual}} hold.
In particular, for both $(\mathrm{s},\mathrm{a})$- and $\mathrm{r}$-rectangular sets, the equation \eqref{equivalence_and_strong_dual} holds.
{\rm (ii)}
For $\mathrm{s}$-rectangular sets, if $\cP_{s_t}^t$ is convex (Assumption {\rm \ref{assump_margin_convex}}) for all $s_t \in \s_t$ and $t = 1, \ldots, T$, then \eqref{equivalence_and_strong_dual} holds.
\end{cor}

\begin{proof}
In view of Theorem \ref{prop_form_eq_and_strong_static_dual} it suffices to verify that a saddle point $(\pi_t^*, P_{t}^*)$ exists for \eqref{mdp-5} and $P_{t}^*$ does not depend on $s_t$.
This trivially holds for ambiguity sets satisfying Assumption \ref{assum}(b),
and in particular, for $(\mathrm{s},\mathrm{a})$- and $\mathrm{r}$-rectangular sets.
On the other hand, for $\mathrm{s}$-rectangular sets satisfying Assumption \ref{assump_margin_convex}, the existence of saddle point $(\pi_t^*, P_{t}^*)$ for \eqref{mdp-5} follows from Theorem \ref{prop_s_rect_convex_marginal}.
The fact that $P_{t}^*$ can be chosen independent of $s_t$ follows trivially using the similar argument in Corollary \ref{prop_sa_rec_strong_dual} and the definition of $\mathrm{s}$-rectangularity.
\end{proof}

Similar to Remark \ref{remark_form_equivalence}, the strong duality of the static formulation involving $(\mathrm{s}, \mathrm{a})$- and $\mathrm{r}$-rectangular sets has been discussed in \cite{Nilim2005} and \cite{goyal2023robust} respectively.
Here our argument differs by first discussing the strong duality of the game formulation and then using formulation equivalence.
Both steps are conceptually simple and provides a unifying treatment of different ambiguity sets.
It is worth noting here that for $\mathrm{s}$-rectangular sets, the strong duality hinges upon the convexity of  $\{\cP_{s_t}^t\}$.
Indeed, it can be readily seen that Example \ref{ex_s_rect_non_saddle} considers an $\mathrm{s}$-rectangular ambiguity set with non-convex $\cP_{s_A}^1$.
Consequently we have
\begin{align*}
 \mathrm{OPT} \eqref{mdp-1}   =  \mathrm{OPT} \eqref{eq_static_formulation}  <   \mathrm{OPT} \eqref{dual-1} =  \mathrm{OPT}   \eqref{eq_static_formulation_dual},
\end{align*}
where the equalities follows from Corollary \ref{cor_form_equivalence}, and the strict inequality follows from Example \ref{ex_s_rect_non_saddle}.
It appears that this essential role of convexity in showing the strong duality of static formulation involving $\mathrm{s}$-rectangular sets is not explicitly discussed in \cite{kuhn2013}.

\begin{remark}
{\rm
Connections between distributionally robust MDPs and stochastic games have been made in several existing works.
In \cite{iyen,Nilim2005}, static formulation with $\mathrm{(s,a)}$-rectangular ambiguity sets is shown to be an instance of perfect information sequential games.
The constructed game therein is asymmetric as the nature chooses its action after the controller.
In \cite{mannor2012lightning}, a non-rectangular ambiguity set termed LDST set is proposed by introducing certain integer-valued budget-constraint (on the deviation to nominal kernel across all stages) of the nature.
Instead of considering coupling among actions or states within the same stage (such as $\mathrm{s}$- or $\mathrm{r}$-rectangularity), therein non-rectangularity mostly refers to the coupling of transition probabilities across different stages (cf., \cite{shaope2021}, Eq. (2.8)) that naturally arises from the budget constraint of the nature. 
Consequently, dynamic equations are defined over an augmented state space including the remaining budget. These equations are established by converting to a perfect information game on the augmented space.
Already for the LDST set, formulation equivalence seems to no longer hold as the static formulation is NP-hard.
A generalization of this non-(time)-rectangular set is made in \cite{mannor2016robust} for finite-horizon MDPs, by considering an asymmetric dynamic game where nature chooses its action after the controller.
Even if the ambiguity set becomes time-rectangular therein, it appears that the connection between game and static formulations remains unclear, and the associated strong duality is difficult to discuss given the asymmetry.

In comparison, the proposed game formulation \eqref{mdp-1} (resp. \eqref{dual-1}) is symmetric in the sense that actions are made simultaneously by the controller and the nature.
We also provide explicit discussions on strong duality in Section \ref{sec-boc} and establish its connections and implications to the static formulation in Section \ref{subsec_static_form}.
}
$\hfill \square$
\end{remark}

\subsection{Construction of   new ambiguity sets.}\label{subsec_sr_rectangular}

We now discuss an application of our analytical approach to the construction of a new class of ambiguity sets.
The proposed ambiguity sets {\it strictly} contain  $\mathrm{(s,a)}$-, $\mathrm{s}$-, and $\mathrm{r}$-rectangular sets.
We also establish formulation equivalence and strong duality for this larger class of ambiguity sets.

\begin{definition}[$\mathrm{sr}$-rectangularity]\label{def_sr_rectangular}
The ambiguity set $\PP_t$ is said to be $\mathrm{sr}$-rectangular if
\begin{align}\label{eq_def_sr_rectangular}
\PP_t = \beta \PP^{\mathrm{s}}_t + (1-\beta) \PP^{\mathrm{r}}_t,
\end{align}
for some $\beta \in [0,1]$,  $\mathrm{s}$-rectangular set $\PP^{\mathrm{s}}_t$,
and  $\mathrm{r}$-rectangular set $\PP^{\mathrm{r}}_t$.
We say $\PP_t$ is $\mathrm{s}$-convex if $\PP^{\mathrm{s}}_t$ has convex state-wise marginalization.
\end{definition}

Conceptually, $\mathrm{sr}$-rectangular sets are richer than $\mathrm{s}$- or $\mathrm{r}$-rectangular sets:
compared to $\mathrm{s}$-rectangularity, $\mathrm{sr}$-rectangularity additionally allows one to consider potential coupling of transition probabilities across different states;
and it allows one to model arbitrary coupling of transition probabilities across different actions compared to $\mathrm{r}$-rectangularity.
We make this argument precise in the following result.

\begin{proposition}
If $\PP_t$ is either $\mathrm{(s,a)}$-, $\mathrm{s}$-, or $\mathrm{r}$-rectangular,
then it is $\mathrm{sr}$-rectangular.
On the other hand, there exists $\PP_t$ that is $\mathrm{sr}$-rectangular but not $\mathrm{(s,a)}$-, $\mathrm{s}$-, or $\mathrm{r}$-rectangular.
\end{proposition}

\begin{figure}[t]
    \centering
    \includegraphics[width=0.68\textwidth]{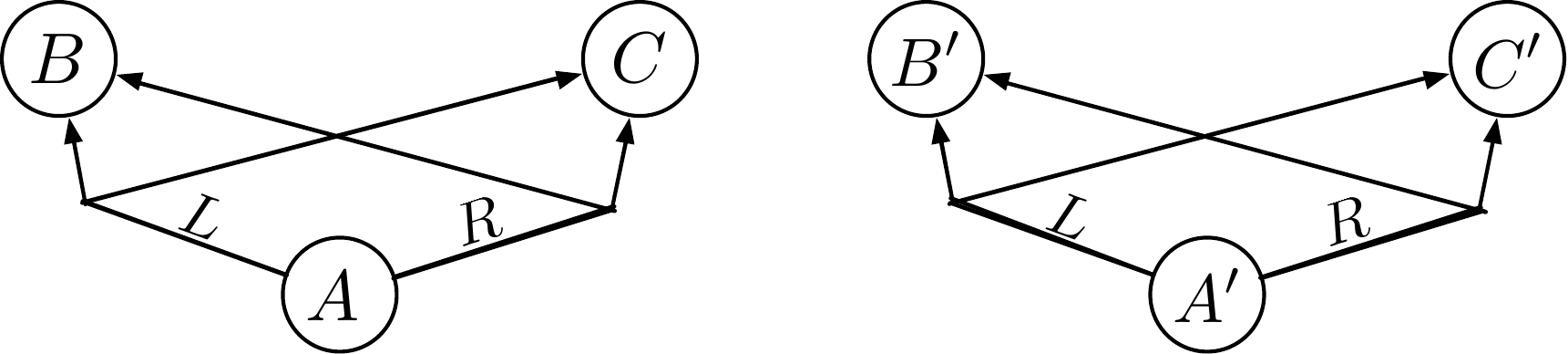} %
    \vspace{-0.1in}
    \caption{Two-stage distributionally robust MDP with $\mathrm{sr}$-rectangular ambiguity set. }
    \label{fig:sr_rectangular}
\end{figure}

\begin{proof}
The first part of the claim trivially holds as any $\mathrm{(s,a)}$-rectangular set is also $\mathrm{s}$-rectangular.
To show the second part of the claim, let us consider the two-stage MDP instance illustrated in Figure \ref{fig:sr_rectangular}.
Similar to Example \ref{ex_example_non_necessary_of_sa}, $\s_1 = \{s_A, s_A' \}$, $\s_2 = \{s_B, s_B', s_C, s_C'\}$  and there  are two actions $a_L$ and $a_R$ at either state $s_A$ or $s_A'$, i.e.   $\A_1=\{a_L,a_R\}$.
 The cost functions are  given by $c_1 \equiv 0$,  $c_2(s_B) = c_2(s_B') =1$ and $c_2(s_C) = c_2(s_C') =0$.
Here $\PP_1 = \frac{1}{2} \PP^{\mathrm{s}}_1+ \frac{1}{2} \PP^{\mathrm{r}}_1$, where
\begin{align*}
\PP^{\mathrm{r}}_1 & = \left\{
P_1:
P_1(s_B|s_A, a_L) = p, P_1(s_B|s_A, a_R) = \tfrac{1+p}{3};
P_1(s_B'|s_A', a_L) = p, P_1(s_B'|s_A', a_R) = \tfrac{1+p}{3};
 p \in [0,1]
\right\}, \\
\PP^{\mathrm{s}}_1 & = \left\{
P_1:
P_1(s_B|s_A, a_L) = q, P_1(s_B|s_A, a_R) = 1-q;
 P_1(s_B'|s_A', a_L) = z, P_1(s_B'|s_A', a_R) = 1-z;  q,z \in [0,1]
\right\}.
\end{align*}
Clearly $\PP^{\mathrm{s}}_1$ is $\mathrm{s}$-rectangular, and $\PP^{\mathrm{r}}_1$ is $\mathrm{r}$-rectangular (see Remark \ref{remark_r_rect}).
We proceed to show that $\PP_1$ is neither $\mathrm{s}$- or $\mathrm{r}$-rectangular.

To show that $\PP_1$ is not $\mathrm{s}$-rectangular.
Suppose one choose $P_1 \in \PP_1$ such that $P_1(s_B|s_A, a_L) = 1$.
This immediately implies $p = q = 1$.
Subsequently we must have $P_1(s_B'|s_A', a_L) \geq \tfrac{1}{2}$.
On the other hand, it is clear that
there exists $\tilde{P}_1$ such that $\tilde{P}_1(s_B' | s_A', a_L) = 0$.
Combining previous two observations implies $\PP_1$ cannot be $\mathrm{s}$-rectangular.

To show that $\PP_1$ is not $\mathrm{r}$-rectangular.
One can verify through direct computation (similar to Example \ref{example_randomized_opt_policy}) that the primal problem \eqref{mdp-1} (with $s_1 = s_A$) has a unique optimal policy given by
\begin{align*}
\pi^*_1(a_L|s_A) = \pi^*_1(a_R|s_A) = 1/2.
\end{align*}
Subsequently $\PP_1$ is not $\mathrm{r}$-rectangular, since otherwise it would contradict with Corollary \ref{prop_r_rectangular} that certifies the existence of a non-randomized optimal policy (through Theorem \ref{th-duality} - (iii)).
\end{proof}

\begin{figure}[t]
    \centering
    \includegraphics[width=0.5\textwidth]{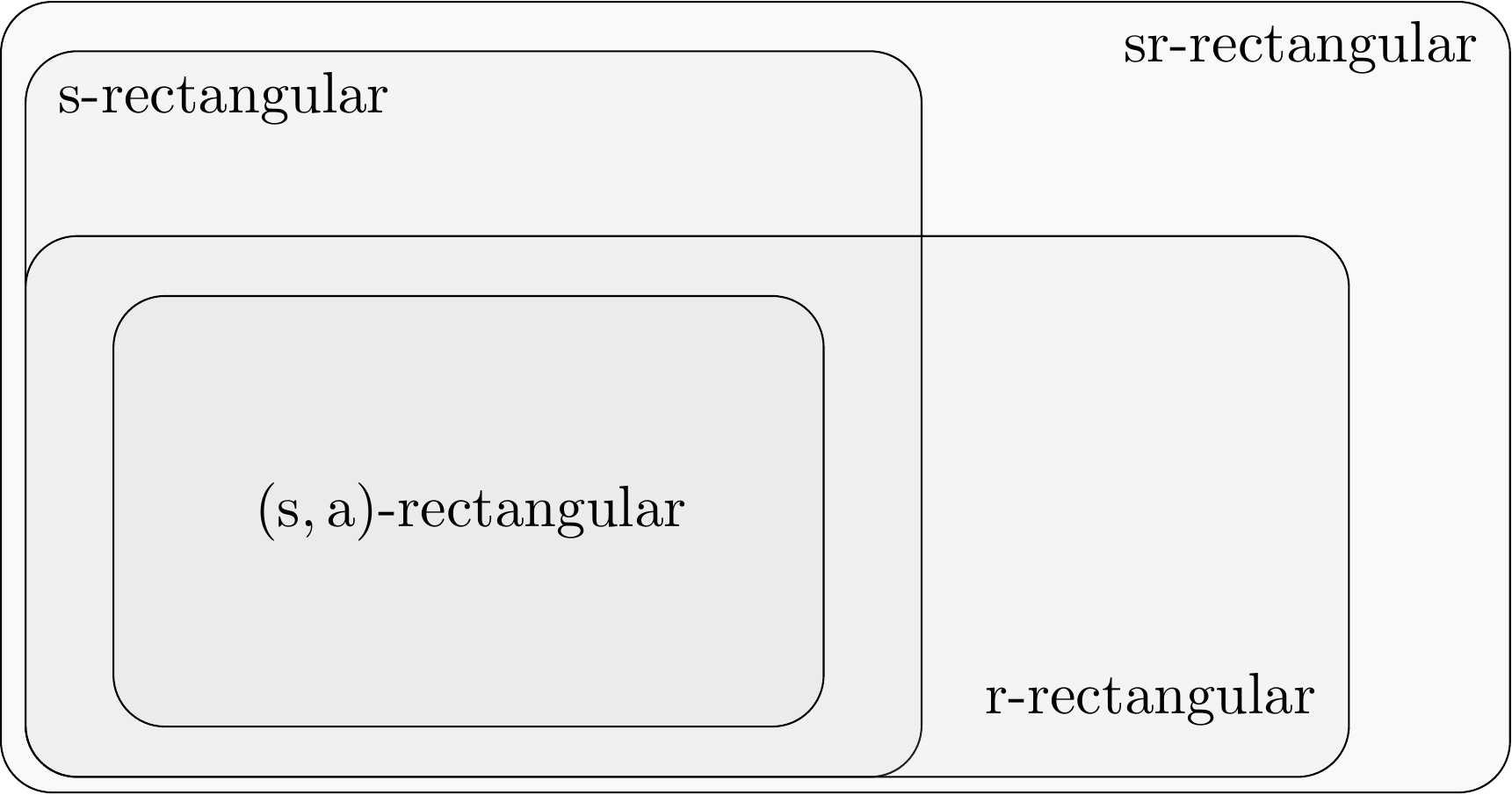} %
    \vspace{-0.1in}
    \caption{Relations between $\mathrm{(s,a)}$-, $\mathrm{s}$-, $\mathrm{r}$-, and $\mathrm{sr}$-rectangular sets.
    It has been shown in \cite{goyal2023robust} that neither the class of $\mathrm{s}$- or $\mathrm{r}$-rectangular sets contains each other. }
    \label{fig:relation}
\end{figure}

Figure \ref{fig:relation} illustrates the proposed $\mathrm{sr}$-rectangularity compared to existing classes of rectangular sets.
We are ready to establish the equivalence of game and static formulations, and the strong duality of two formulations for $\mathrm{sr}$-rectangular sets.

\begin{cor}
For  $\mathrm{sr}$-rectangular sets, the game formulation is equivalent to the static formulation. In particular, \eqref{mdp-5} provides the dynamic equations of the primal problem \eqref{eq_static_formulation} of the static formulation,
and \eqref{dual-4} provides the dynamic equations of the dual problem \eqref{eq_static_formulation_dual} of the static formulation.

In addition, if $\PP_t$ is $\mathrm{s}$-convex, then strong duality for both the game and static formulations holds (i.e.,  \eqref{equivalence_and_strong_dual} holds),
and in this case there exists an optimal policy of the controller in \eqref{prob_hist_controller_markov_nature} that is Markovian.
\end{cor}

\begin{proof}
Proof follows from similar argument as in Corollary  \ref{cor_existence_markov_opt_policy},  \ref{cor_form_equivalence}, and \ref{cor-equality}.
\end{proof}

We conclude this section by noting that the $\mathrm{sr}$-rectangularity generalizes existing notions of rectangularity \cite{iyen,Nilim2005,goyal2023robust, goh2018data, Tallec, kuhn2013} for the static formulation while admitting natural dynamic equations.
Such dynamic equations arise, again, precisely due to the equivalence between static and game formulations under the specified rectangularity of the associated ambiguity sets.

\setcounter{equation}{0}
\section{Cost-robust MDP}\label{sec_cost_robust}
Our discussion so far has focused on the ambiguity of the transition kernels.
In this section, we show that by following virtually the same lines of argument, one can naturally consider ambiguity on the cost functions.
Specifically,   we assume in this section that  the transition kernels $\{P_t\}$ are given (fixed), and  there exists a set $\C_t$, $t = 1, \ldots, T$, of cost functions.
For simplicity we assume that there is no cost ambiguity at $t = T+1$.
It is worth noting   that since each $c_t \in \C_t$ is of the same dimension as the $P_t$,
we can define $(\mathrm{s,a})$- and $\mathrm{s}$-rectangularity for $\C_t$ in exactly the same manner as in Sections \ref{sec-boc} and \ref{subsec_static_form}.

\subsection{Game formulation of cost-robust MDP}
We define the game formulation of cost-robust MDPs similarly as in Section \ref{sec-boc}.
In this case, a (Markovian) policy of the nature maps from $\s_t$ into $\C_t$. As before
we denote by $\Pi$ and $\Gamma$ the sets   of Markovian policies of the controller and the nature, respectively.
The corresponding primal problem of the game formulation is
\begin{equation}\label{mdp-1_cost}
 \min_{\pi\in \Pi}\max_{\gamma\in \Gamma}\bbe^{\pi,\gamma}
 \left [ \tsum_{t=1}^{T}
c_t(s_t,a_t,s_{t+1})+c_{T+1}(s_{T+1})
\right].
\end{equation}
Given a policy $\{\pi_t\}$ of the controller, the nature chooses a policy
so as to maximize the total expected cost, i.e.,
\begin{equation}\label{mdp-2_cost}
\max_{\gamma \in \Gamma}  \bbe^{\pi, \gamma}\left [ \tsum_{t=1}^{T}
c_t(s_t,a_t,s_{t+1})+c_{T+1}(s_{T+1})
\right].
\end{equation}

The following propositions follow the same argument as that of Proposition \ref{pr-dyncontr} and \ref{prop_game_primal_dp}.

\begin{proposition}
\label{pr-dyncontr_cost}
Given a policy $\pi\in \Pi$ of the controller, problem
 \eqref{mdp-2_cost} admits the following
 dynamic programming equations: $\V^\pi_{T+1}(s_{T+1})=   c_{T+1}(s_{T+1})$ and for $t=T,...,1$,
 \begin{align}
 \label{mdp-3_cost}
\V^\pi_t(s_t)  & = \max_{c_t \in \C_t}  \tsum_{s_{t+1} \in \s_{t+1}}  \tsum_{a_t \in \A_t(s_t)} P_t(s_{t+1} | s_t ,a_t) \pi_t(a_t | s_t) \big[
c_t(s_t, a_t ,s_{t+1})+\V^\pi_{t+1}(s_{t+1}) \big],
\end{align}
where
 $\V^\pi_t(s_t)$  corresponds to the optimal value of problem
\begin{align*}
\max_{\gamma \in \Gamma}  \bbe^{\pi, \gamma}\left [ \tsum_{i=t}^{T}
c_i(s_i,a_i,s_{i+1})+c_{T+1}(s_{T+1})
\right].
\end{align*}
Furthermore, it holds that
\begin{align}
\label{cost_dp_in_regularization}
\V^\pi_t(s_t)
& = h_{s_t}(\pi_t(\cdot|s_t))
+ \tsum_{s_{t+1} \in \s_{t+1}}  \tsum_{a_t \in \A_t(s_t)}   P_t(s_{t+1} | s_t ,a_t) \pi_t(a_t|s_t) \V^\pi_{t+1}(s_{t+1}),
\end{align}
where
\begin{equation}
\label{h-func}
h_{s_t}(x_t) := \sup_{{c}_t \in {\C}_t}    \tsum\limits_{a_t \in \A_t(s_t)} x_t(a_t) \tsum\limits_{s_{t+1} \in \s_{t+1}} P_t(s_{t+1} | s_t ,a_t)  c_t(s_t, a_t, s_{t+1}).
\end{equation}
\end{proposition}

It can be noted that function $h_{s_t}(\cdot)$, defined   in \eqref{h-func}, is given by the maximum of  linear functions in $x_t\in \bbr^{|\A_t(s_t)|}$, and hence is convex and positively homogeneous for any $s_t\in \s_t$.

\begin{proposition}\label{prop_game_primal_dp_cost}
The dynamic programming equations for the primal  problem \eqref{mdp-1_cost} are: $V_{T+1}(s_{T+1})=  c_{T+1}(s_{T+1})$ and for $t=T,...,1$,
\begin{equation}
\label{mdp-5_cost}
V_t(s_t)= \min_{\pi_t(\cdot|s_t) \in \Delta_{|\A_t(s_t)|}} \max_{c_t \in \C_t} \tsum\limits_{s_{t+1} \in \s_{t+1}} \tsum\limits_{a_t \in \A_t(s_t)} P_t(s_{t+1} | s_t ,a_t) \pi_t(a_t | s_t) \big[
c_t(s_t, a_t ,s_{t+1})+V_{t+1}(s_{t+1}) \big],
\end{equation}
where $V_t(s_t)$ corresponds to the optimal value of
 \begin{align*}
 \min_{\pi \in \Pi}
 \max_{\gamma \in \Gamma}  \bbe^{\pi, \gamma}\left [ \tsum_{l=t}^{T}
c_l(s_l,a_l,s_{l+1})+c_{T+1}(s_{T+1})
\right].
 \end{align*}
 Furthermore, it holds that
\begin{align}
\label{cost_dp_in_regularization-2}
V_t(s_t)
& = \min_{\pi_t(\cdot|s_t) \in \Delta_{|\A_t(s_t)|}}  \left\{ h_{s_t}(\pi_t(\cdot|s_t))
+ \tsum_{s_{t+1} \in \s_{t+1}}  \tsum_{a_t \in \A_t(s_t)}   P_t(s_{t+1} | s_t ,a_t) \pi_t(a_t|s_t) V_{t+1}(s_{t+1})\right\}.
\end{align}
\end{proposition}

It follows  from  \eqref{cost_dp_in_regularization-2} in Proposition \ref{prop_game_primal_dp_cost}
 that the game formulation of any cost-robust MDP can be reduced to a regular (i.e., non-robust) MDP with  convex policy-dependent cost functions, provided the function $h_{s_t}(\cdot)$ is real valued for every $s_t\in\s_t$.

 \begin{theorem}
 The game formulation \eqref{mdp-1_cost} of any cost-robust MDP is equivalent to a regular (i.e., non-robust) MDP with the same state and action spaces and transition kernels.
For a given policy $\pi \in \Pi$, the cost of the regular MDP at state $s_t \in \s_t$ is given by
$h_{s_t}(\pi_t(\cdot|s_t))$.
 \end{theorem}

 The dual of problem \eqref{mdp-1_cost} is defined as
\begin{equation}\label{dual-1_cost}
 \max_{\gamma\in \Gamma} \min_{\pi\in \Pi} \bbe^{\pi,\gamma}
 \left [ \tsum_{t=1}^{T}
c_t(s_t,a_t,s_{t+1})+c_{T+1}(s_{T+1})
\right].
\end{equation}
The  dynamic programming equations for the dual  problem are: $Q_{T+1}(s_{T+1})=  c_{T+1}(s_{T+1})$ and for $t=T,...,1$,
  \begin{equation}\label{dual-4_cost}
Q_t(s_t)=\max_{c_t\in \C_t} \min_{\pi_t(\cdot|s_t) \in \Delta_{|\A_t(s_t)|}}   \tsum\limits_{s_{t+1} \in \s_{t+1}} \tsum\limits_{a_t \in \A_t(s_t)} P_t(s_{t+1} | s_t ,a_t) \pi_t(a_t | s_t) \big[
c_t(s_t, a_t ,s_{t+1})+Q_{t+1}(s_{t+1}) \big].
\end{equation}

Similar to Section \ref{sec-boc}, we now introduce sufficient conditions that certify the strong duality of the game formulation for cost-robust MDP.

 \begin{assumption}
\label{assum_cost}
\begin{enumerate}
\item [{\rm (a)}]
 For every $s_t\in \s_t$ there is a cost function $c^*_{t}\in \C_t$, potentially depending on $s_t$, with
\begin{equation}\label{saddle-b_cost}
c^*_{t} \in
\argmax\limits_{c_t \in \C_t} \tsum_{s_{t+1} \in \s_{t+1}} \tsum_{a_t \in \A(s_t)} P_t(s_{t+1} | s_t ,a_t) \mu_t(a_t) \big[
c_t(s_t, a_t ,s_{t+1})+V_{t+1}(s_{t+1}) \big],
\end{equation}
for   any $\mu_t \in \Delta_{|A_t(s_t)|}$.
\item [{\rm (b)}]
There is a cost function $c^*_{t}\in \C_t$
  such that condition \eqref{saddle-b_cost}  holds for all $s_t\in \s_t$  and any $\mu_t \in \Delta_{|A_t(s_t)|}$.
   \end{enumerate}
\end{assumption}

Suppose that  Assumption \ref{assum_cost}(a) holds   and consider
\begin{equation}\label{at_cost}
\pi_t^*(\cdot|s_t)=\delta_{a_t^*(s_t)} \; \text{ with}\;a_t^*(s_t)\in\argmin_{a_t\in \A_t(s_t)} \tsum\limits_{s_{t+1} \in \s_{t+1}} P_{t}(s_{t+1} | s_t ,a_t)   \big[
c^*_t(s_t, a_t ,s_{t+1})+V_{t+1}(s_{t+1}) \big].
\end{equation}
   It follows directly from \eqref{saddle-b_cost} and  \eqref{at_cost}
that    $(\pi_t^* ,c_{t}^*)$ is a saddle point of the min-max problem  \eqref{mdp-5_cost}.
Consequently  we have the following result.

\begin{theorem}
\label{th-duality_cost}
Suppose that Assumption {\rm   \ref{assum_cost}(a)} is fulfilled for all $t=1,...,T$,  and let $\pi_t^*(\cdot|s_t)$ be a minimizer defined in \eqref{at_cost}. Then the following holds:
\begin{itemize}[itemsep=0.5pt, topsep=3pt]
\item[{\rm (i)}] $(\pi_t^*,c_{t}^*)$ is a saddle point of the min-max problem  \eqref{mdp-5_cost}.
\item[{\rm (ii)}] There is no duality gap between  the primal problem \eqref{mdp-1_cost} and its dual  \eqref{dual-1_cost}.
\item[{\rm (iii)}]
$\{\pi^*_t\}$ is an optimal (non-randomized)  policy of the controller in the primal problem \eqref{mdp-1_cost} and
$
 \gamma_t^*(s_t) = c_{t}^*$,  $ t = 1, \ldots, T,
 $
 is an optimal policy of the nature considered in the dual problem \eqref{dual-1_cost}  (here $c_t^*$, viewed as an element of $\C_t$,   can depend on $s_t\in \s_t$).
 \item[{\rm (iv)}]  If moreover  Assumption {\rm   \ref{assum_cost}(b)} holds, then the above claims hold with $c_{t}^*$ independent of $s_t$.
 \end{itemize}
\end{theorem}

An immediate application of Theorem \ref{th-duality_cost} yields the following propositions.
\begin{cor}\label{prop_sa_rec_strong_dual_cost}
Suppose that the ambiguity sets $\C_t$, $t=1,...,T$, are
$(\mathrm{s}, \mathrm{a})$-rectangular.
Then Assumption {\rm  \ref{assum_cost}(b)} holds,
and conclusions {\rm (i) - (iv)} of Theorem {\rm \ref{th-duality_cost}} follow.
\end{cor}

The following condition can be viewed as the counterpart of Assumption \ref{assump_margin_convex} for the cost functions.

\begin{assumption}\label{assump_margin_convex_cost}
The state-wise marginalization of $\C_t$, defined as
\begin{align}\label{eq_marginal_s_cost}
\cC^t_{s_t}: = \{
q_t: q_t (\cdot,\cdot)= c_t(\cdot|s_t, \cdot) ~ \text{for some} ~ c_t \in \C_t
\},
\end{align}
is convex.
\end{assumption}

 Clearly, Assumption \ref{assump_margin_convex_cost} holds if the set $\C_t$ is convex.
The next result then follows from the same argument as in Theorem \ref{prop_s_rect_convex_marginal}.

\begin{theorem}\label{prop_s_rect_convex_marginal_cost}
Suppose Assumption {\rm  \ref{assump_margin_convex_cost}} holds.
Then
\begin{itemize}[itemsep=0.5pt, topsep=3pt]
\item[{\rm (i)}] There exists a saddle point $(\pi_t^*,c_{t}^*)$ for the min-max problem  \eqref{mdp-5_cost}.
\item[{\rm (ii)}] There is no duality gap between  the primal problem \eqref{mdp-1_cost} and its dual  \eqref{dual-1_cost}.
 \end{itemize}
Furthermore, the primal problem \eqref{mdp-1_cost} has a non-randomized optimal policy iff the min-max problem
\begin{equation*}
  \min_{a_t \in  \A_t(s_t)} \max_{c_t \in \C_t} \tsum\limits_{s_{t+1} \in \s_{t+1}} P_t(s_{t+1} | s_t ,a_t)\big[
c_t(s_t, a_t ,s_{t+1})+V_{t+1}(s_{t+1}) \big].
\end{equation*}
 has a saddle point for all $t=1,...,T$ and $s_t\in \s_t$.
\end{theorem}

\subsection{Static formulation of cost-robust MDP}
Let us consider the following  static formulation of cost-robust MDP,
\begin{align}\label{eq_static_formulation_cost}
 \min_{\pi\in \Pi}\sup_{c_1\in \C_1, \ldots, c_T \in \C_T} \bbe^{\pi,\{P_t\}_{t = 1}^T}
 \left [ \tsum_{t=1}^{T}
c_t(s_t,a_t,s_{t+1})+c_{T+1}(s_{T+1})
\right].
\end{align}
The corresponding dual problem is given by
\begin{align}\label{eq_static_formulation_dual_cost}
\max_{c_1\in \C_1, \ldots, c_T \in \C_T} \inf_{\pi\in \Pi} \bbe^{\pi,\{P_t\}_{t = 1}^T}
 \left [ \tsum_{t=1}^{T}
c_t(s_t,a_t,s_{t+1})+c_{T+1}(s_{T+1})
\right].
\end{align}
Formulations \eqref{eq_static_formulation_cost} and \eqref{eq_static_formulation_dual_cost} can be viewed as
 the static counterparts of the respective problems  \eqref{mdp-1_cost} and \eqref{dual-1_cost}.

 With the same argument as in  Theorem \ref{prop_form_equivalence}, one can establish the  following equivalence of game and static formulations for cost-robust MDP.

\begin{theorem}\label{prop_form_equivalence_cost}
The following holds.
{\rm (i)}  If for any policy $\pi\in \Pi$,  there exists a solution $c_{t}^*$ of \eqref{mdp-3_cost} that is independent of
 $s_t$, for every $s_t \in \s_t$ and $t = 1,\ldots, T$,
then the primal problem \eqref{mdp-1_cost} of the game  and its static counterpart    \eqref{eq_static_formulation_cost}   are equivalent, i.e.,
\begin{align}
 \mathrm{OPT} \eqref{mdp-1_cost} =  \mathrm{OPT} \eqref{eq_static_formulation_cost},
  \end{align}
  and the dynamic equations \eqref{mdp-5_cost} hold for the static problem  \eqref{eq_static_formulation_cost} as well.

{\rm (ii)}  If for problem  \eqref{dual-4_cost} (not necessarily having a saddle point), there exits a solution $c_{t}^*$ such that $c_{t}^*$ is independent of
 $s_t$, for every $s_t \in \s_t$ and $t = 1,\ldots, T$,
then the dual problem   \eqref{dual-1_cost} of the game and its  static counterpart   \eqref{eq_static_formulation_dual_cost}    are equivalent, i.e.,
\begin{align}
   \mathrm{OPT} \eqref{dual-1_cost} = \mathrm{OPT} \eqref{eq_static_formulation_dual_cost} ,
\end{align}
and the dynamic equations \eqref{dual-4_cost} hold for the static problem  \eqref{eq_static_formulation_dual_cost} as well.
\end{theorem}

\begin{cor}\label{cor_form_equivalence_cost}
For  the $(\mathrm{s},\mathrm{a})$- and $\mathrm{s}$-rectangular sets, the game formulation is equivalent to the static formulation for cost-robust MDP. In particular, \eqref{mdp-5_cost} provides the dynamic equations of the primal problem \eqref{eq_static_formulation_cost} of the static formulation,
and \eqref{dual-4_cost} provides the dynamic equations of the dual problem \eqref{eq_static_formulation_dual_cost} of the static formulation.
\end{cor}

We now conclude this section by discussing the duality of the static formulation.
As before, the following result follows from the same lines as that of Section \ref{sec-boc} (i.e., Theorem \ref{prop_form_eq_and_strong_static_dual}).

\begin{theorem}\label{prop_form_eq_and_strong_static_dual_cost}
Suppose for every $t = 1,\ldots T$ and $s_t \in \s_t$, a saddle point $(\pi_t^*, c_{t}^*)$ exists for \eqref{mdp-5_cost}, and $c_{t}^*$ does not depend on $s_t$, then
\begin{align}\label{equivalence_and_strong_dual_cost}
 \mathrm{OPT} \eqref{mdp-1_cost} =  \mathrm{OPT} \eqref{dual-1_cost}  =  \mathrm{OPT} \eqref{eq_static_formulation_cost} = \mathrm{OPT} \eqref{eq_static_formulation_dual_cost}.
\end{align}
That is,  the static formulation (for both primal and dual problems) is equivalent to the game formulation,
and the strong duality holds for the static formulation.
\end{theorem}

\begin{cor} {\rm (i)}
Under Assumption {\rm \ref{assum_cost}(b)}, the conclusions of Theorem {\rm \ref{prop_form_eq_and_strong_static_dual_cost}} hold.
In particular, for $(\mathrm{s},\mathrm{a})$-rectangular sets, \eqref{equivalence_and_strong_dual_cost} holds.
{\rm (ii)}
For $\mathrm{s}$-rectangular sets, if $\cC_{s_t}^t$ is convex (Assumption {\rm \ref{assump_margin_convex_cost}}) for all $s_t \in \s_t$ and $t = 1, \ldots, T$, then \eqref{equivalence_and_strong_dual_cost} holds.
\end{cor}

It is worth noting that for $\mathrm{s}$-rectangular ambiguity sets, equivalence between cost-robust MDPs and regularized MDPs has also been discussed in \cite{derman2021twice}.
Our analytical approach here in addition provides a natural discussion on the strong duality of cost-robust MDPs via a unified treatment similar to distributionally robust MDPs discussed in Sections  \ref{sec-boc} and \ref{subsec_static_form}.

\setcounter{equation}{0}
\section{Risk averse setting}\label{sec_risk_averse}

In this section we discuss a risk averse formulation of  MDPs and its relation to the  distributionally robust approach. Consider a finite set $\O=\{\w_1,...,\w_n\}$ equipped with sigma algebra $\F$ of all subsets of $\O$. Let $\Z$ be the space of functions $Z:\O\to \bbr$. Note that $Z$ can be identified with $n$-dimensional  vector with components $Z_i=Z(\w_i)$, $i=1,...,n$, and hence $\Z$ can be identified with $\bbr^n$. Let $\R:\Z\to \bbr$ be a convex  functional. Suppose that $\R$ is monotone, i.e.  if $Z\le Z'$, then $\R(Z)\le \R(Z')$; positively homogeneous, i.e. $\R(\lambda Z)=\lambda \R(Z)$ for any $\lambda\ge 0$ and $Z\in \Z$; and such that $\R(Z+a)=\R(Z)+a$ for any $a\in \bbr$ and $Z\in \Z$. The functional $\R$ can be viewed as a coherent risk measure defined on the space $\Z$ of (measurable) functions (cf., \cite{ADEH:1999}). We can refer to \cite{fol04,SDR}  for a thorough discussion of coherent risk measures defined on general probability spaces.

It is said that $\R=\R_P$  is law   invariant with respect to a probability measure $P$,  defined on $(\O,\F)$, if $\R_P(Z)$ can be considered as a function of the cumulative distribution function $F_Z(z)=P(Z\le z)$ (e.g., \cite[Section 6.3.3]{SDR}). Law invariant coherent risk measure $\R_P:\Z\to \bbr$ has the following dual representation,  in the distributionally robust form,
\begin{equation}
\label{mesdual}
\R_P (Z) = \sup_{q\in \cM_P}\tsum_{i=1}^n  q_i Z_i,
\end{equation}
where $\cM_P\subset \Delta_n$ is a set of probability measures (probability vectors), depending on $P$, such that  every $Q\in \cM_P$ is absolutely continuous with respect to $P$. Note that  in the considered setting of finite $\O$, measure  $Q$ is absolutely continuous with respect to $P$ iff $p_i=0$ implies that $q_i=0$.  Of course if all probabilities  $p_i$  are positive, any probability
measure $Q$  is absolutely continuous with respect to $P$.

Let $\bbp_t$ be a specified (reference)  transition kernel. Considering  $\bbp_t(\cdot|s_t,a_t)$  as a probability measure on $\s_{t+1}$,  we can consider the corresponding law invariant coherent  risk measure  $\R_{\bbp_t(\cdot|s_t,a_t)}$ depending on $s_t\in \s_t$ and $a_t\in \A_t(s_t)$.
 This suggests  the following counterpart of dynamic equations \eqref{mdp-5}:
\begin{equation}\label{risk-1a}
V_t(s_t)= \min_{\pi_t(\cdot|s_t) \in \Delta_{|\A(s_t)|}} \tsum_{a_t \in \A_t(s_t)}  \pi_t(a_t | s_t)   \R_{\bbp_t(\cdot|s_t, a_t)} \big(
c_t(s_t, a_t ,\cdot)+V_{t+1}(\cdot) \big).
\end{equation}
From the dual representation of $\R_{\bbp_t(\cdot|s_t,a_t)}$, the dynamic equation \eqref{risk-1a} can be written in the form of \eqref{mdp-5}:
\begin{equation}
\label{risk-dyn}
\begin{array}{lll}
V_t(s_t) & =& \min\limits_{\pi_t(\cdot|s_t) \in \Delta_{|\A_t(s_t)|}} \sup\limits_{\;q \in \cM_{\bbp_t(\cdot|s_t, a_t)}} \tsum\limits_{s_{t+1} \in \s_{t+1}}  \tsum\limits_{a_t \in \A_t(s_t)} q(s_{t+1})\pi_t(a_t | s_t) [
c_t(s_t,a_t,s_{t+1})+V_{t+1}(s_{t+1})]
 \\
& = & \min\limits_{\pi_t(\cdot|s_t) \in \Delta_{|\A_t(s_t)|}} \sup\limits_{\; P_t\in \PP_t} \tsum\limits_{s_{t+1} \in \s_{t+1}}  \tsum\limits_{a_t \in \A_t(s_t)} P_t(s_{t+1}|s_t, a_t)\pi_t(a_t | s_t) [
c_t(s_t,a_t,s_{t+1})+V_{t+1}(s_{t+1})],
\end{array}
\end{equation}
where
\begin{equation}
\label{risk amb}
\PP_t := \left\{ P_t: P_t(\cdot|s_t, a_t) \in  \cM_{\bbp_t(\cdot|s_t, a_t)}, ~ \forall s_t \in \s_t, ~ \forall a_t \in \A_t(s_t)
\right\}.
\end{equation}
 In a general abstract form such approach to construction of risk averse MDPs was introduced  in Ruszczy\'nski \cite{Rusz2010}.

\begin{remark}
\label{rem-stat}
{\rm
It follows directly from Definition \ref{def-rect} that the  ambiguity sets $\PP_t$ are
 $(\mathrm{s}, \mathrm{a})$-rectangular. Consequently by Theorem \ref{prop_form_equivalence} we have that the corresponding problem of the form \eqref{mdp-1} is equivalent to its static counterpart \eqref{eq_static_formulation} with the respective ambiguity sets of the form  \eqref{risk amb} and dynamic equations \eqref{risk-1a}. Also we have by Theorem \ref{th-duality}(iii) that it suffices here to consider non-randomized policies $a_t=\pi_t(s_t)$, $t=1,...,T$.
 } $\hfill \square$
 \end{remark}

 An important  example of law invariant coherent  risk measure is   the Average Value-at-Risk (also called Conditional Value-at-Risk, Expected Shortfall, Expected Tail Loss):
\[
\avr_{\alpha,P} (Z)=\inf_{\tau\in \bbr}\{\tau+\alpha^{-1}\bbe_P[Z-\tau]_+\},\;\alpha\in (0,1].
\]
The $\avr_{\alpha,P}$ has dual representation \eqref{mesdual} with
$$
\cM_P:=\left\{q\in \bbr^n: 0\le q_i\le    \alpha^{-1} p_i,\;i=1,...,n,\;\tsum_{i=1}^n  q_i=1\right\},
$$
and hence
\begin{equation}\label{risk-2}
\PP_t=
\left\{P_t:  0 \leq P_t(\cdot|s_t,a_t)\le \alpha^{-1}\bbp_t(\cdot|s_t,a_t), ~ \forall s_t \in \s_t, a_t \in \A_t(s_t)
\right\}.
\end{equation}
This leads to  the respective  dynamic programming equations of the form \eqref{risk-dyn} (cf., \cite{dingfein}).
As it is pointed in Remark \ref{rem-stat} the above nested risk averse formulation is  equivalent to its
static counterpart.  Note however that this static formulation is not equivalent to minimizing $\avr_\alpha$ of the total sum of the costs.

 Now let us consider the following construction.
Let $a_t=\pi_t(s_t)$, $t=1,...,T$,  be a non-randomized policy, and
 consider the scenario tree formed by scenarios (sample paths) $s_1,...,s_{T+1}$  of moving from state $s_t\in\s_t$ to state $s_{t+1}\in \s_{t+1}$,  $t=1,...,T$. The corresponding   probability distribution on this scenario tree is  determined by the conditional probabilities $\bbp_t(s_{t+1}|s_t,\pi_t(s_t))$. In turn this defines the {\em  nested} functional $\cR^\pi:\Z_1\times\cdots\times \Z_{T+1}\to \bbr$ associated with the coherent risk measure   and the (non-randomized) policy $\pi\in \Pi$, defined on the scenario tree of histories of the states (cf., \cite{Rusz2010}).
  That is, at time period $t$ the conditional risk functional is determined  by the considered risk measure and the reference distribution
  $\bbp_t(\cdot|s_t,\pi_t(s_t))$. Then $\cR^\pi$ is defined as the composition of these conditional functionals  (cf.,\cite{shaope2021}, \cite[Section 6.5.1, equations (6.225) and (6.227)]{SDR}).
  Consequently in the considered risk averse setting the  (primal)  distributionally robust  problem
   can be written as
 \begin{equation}\label{risk-1}
 \min_{\pi\in \Pi} \cR^{\pi}
 \left [ \tsum_{t=1}^{T}
c_t(s_t,a_t,s_{t+1})+c_{T+1}(s_{T+1})
\right].
\end{equation}
where $\Pi$ consists of non-randomized Markovian policies.
 Of course if $\R_P=\bbe_P$ (the corresponding  set $\cM_P=\{P\}$ is the singleton), then \eqref{risk-1} becomes the risk neutral problem
 \eqref{intra-1}.
 In the  $(\mathrm{s}, \mathrm{a})$-rectangular case, in particular in the above risk averse setting,   the nested formulation \eqref{risk-1} is equivalent to the static formulation \eqref{eq_static_formulation} with the ambiguity sets $\PP_t$ of the form  \eqref{risk amb}.

 \begin{remark}
 \label{rem-risk}
 {\rm
 Note that the probability structure on the scenario tree determined by the conditional distributions $\bbp_t(\cdot|s_t,\pi_t(s_t))$ is Markovian, i.e., the probability $\bbp_t(s_{t+1}|s_t,\pi_t(s_t))$ does not depend on the history $s_1,...,s_{t-1}$. Therefore  if in the nested formulation \eqref{risk-1} the policies $\pi$ are allowed to be history-dependent,  problem \eqref{risk-1}  still possesses  a Markovian optimal policy (compare with Proposition  \ref{prop_hist}(ii)).
}  $\hfill \square$
 \end{remark}

\setcounter{equation}{0}
\section{Optimal control model}
\label{sec-soc}

The MDP setting can be compared with  the   Stochastic Optimal Control  (SOC)  (discrete time, finite horizon) model (e.g., \cite{ber78}):
\begin{equation}\label{soc-1}
\min\limits_{\pi\in \Pi}  \bbe \left [ \tsum_{t=1}^{T}
c_t(s_t,a_t,\xi_t)+c_{T+1}(s_{T+1})
\right].\\
\end{equation}
 Here $\xi_1,...,\xi_T$,   is a sequence of independently distributed random vectors,  the state spaces $\s_t$ and  control (action) spaces $\A_t$ are finite,   the support $\Xi_t$ of distribution    of $\xi_t$ is finite,  and $\Pi$ is a set of randomized policies determined by the state equations
   \begin{equation}\label{state}
s_{t+1}=F_t(s_t,a_t,\xi_t),\;\;t=1,...,T.
   \end{equation}
   We   assume
   that the  distribution    of $\xi_t$   {\em does not depend} on our actions,
$t=1,...,T$. Also for the sake of simplicity   assume that for all $t$  the  action spaces $\A_t$ do not depend on $s_t$, and that the following feasibility condition holds
\[
\big \{s_{t+1}: s_{t+1}=F_t(s_t,a_t,\xi_t),\;s_t\in \s_t,\;a_t\in \A_t,\;\xi_t\in \Xi_t\big\}\subseteq \s_{t+1},\;t=1,...,T.
\]
The SOC setting    can be formulated in the MDP framework by defining the distribution of transition kernel $P^{Q_t}_t(\cdot|s_t,a_t)$
  as
\begin{equation}\label{index}
 P^{Q_t}_t(s_{t+1}\in A|s_t,a_t)=Q_t\{\xi_t\in \Xi_t:F_t(s_t,a_t,\xi_t)\in A\}, \;A\subseteq \s_{t+1},
\end{equation}
where
$Q_t$ is the probability measure of $\xi_t$, $t=1,...,T$.

The distributionally robust   counterpart of \eqref{soc-1} is obtained by considering a set $\cQ_t$ of   probability measures (distributions) of $\xi_t$, supported on $\Xi_t$, $t=1,...,T$.
That is, the ambiguity set of probability distributions of $(\xi_1,...,\xi_T)$ consists of $Q_1\times\cdots\times Q_T$, with $Q_t\in \cQ_t$, $t=1,...,T$. This defines the corresponding ambiguity sets $\PP_t$ of the respective MDP formulation by considering sets of transition kernels determined by equation \eqref{index} for $Q_t\in \cQ_t$.
The corresponding dynamic programming equations     take the form: $V_{T+1}(s_{T+1})=c_{T+1}(s_{T+1})$, and for $t=T,...,1,$
  \begin{equation}\label{contd-2}
V_t(s_t)=\min_{\pi_t(\cdot|s_t)\in \Delta_{|A_t|}}\sup_{Q_t\in \cQ_t}\underbrace{\tsum_{a_t\in \A_t}\tsum_{\xi_t\in \Xi_t}\pi_t(a_t|s_t) {Q(\xi_t)}[c_t(s_t,a_t,\xi_t)+V_{t+1}(F_t(s_t,a_t,\xi_t))]}_
{\bbe^{\pi_t,Q_t}
\left[ c_t(s_t,a_t,\xi_t)+V_{t+1}(F_t(s_t,a_t,\xi_t)) \right]}.
 \end{equation}
Without loss of generality we can make the following assumption.
\begin{assumption}
The sets $\cQ_t$, $t=1,...,T$,  are {\em convex} and {\em closed}.
\end{assumption}

The dual of the min-max problem  \eqref{contd-2} is
 \begin{equation}\label{contd-3}
\max_{Q_t\in \cQ_t} \inf_{\pi_t(\cdot|s_t)\in \Delta_{|A_t|}} \tsum_{a_t\in \A_t}\tsum_{\xi_t\in \Xi_t}\pi_t(a_t|s_t)q_{\xi_t}[c_t(s_t,a_t,\xi_t)+V_{t+1}(F_t(s_t,a_t,\xi_t))].
 \end{equation}
By Sion's theorem   problems
\eqref{contd-2} and \eqref{contd-3} have the same optimal value and possess saddle point $(\pi_t^*,Q^*_t)$.
Note that the minimum in \eqref{contd-3} is attained at a single point of $\Delta_{|A_t|}$, i.e., problem \eqref{contd-3} can be written as
 \begin{equation}\label{contd-6}
  \max_{Q_t\in \cQ_t}\min_{a_t\in  \A_t}  \bbe^{Q_t} [c_t(s_t,a_t,\xi_t)+V_{t+1}(F_t(s_t,a_t,\xi_t))].
 \end{equation}
  Therefore we have the following condition for existence of non-randomized policies of the  distributionally robust SOC problem  associated with the dynamic equations \eqref{contd-2}.

  \begin{proposition}
  The distributionally robust SOC problem possesses a non-randomized optimal  policy iff for every $t=1,...,T$,  and $s_t\in \s_t$ the min-max problem \eqref{contd-6} has a saddle point $(a^*_t,Q^*_t)\in (\A_t,\cQ_t)$.
 \end{proposition}

For the non-randomized policies the dynamic equations \eqref{contd-2}  become
$V_{T+1}(s_{T+1})=c_{T+1}(s_{T+1})$, and for $t=T,...,1,$
 \begin{equation}\label{contdyn}
V_t(s_t)=\inf_{a_t\in \A_t}\sup_{Q_t\in \cQ_t}\bbe^{Q_t}
[c_t(s_t,a_t,\xi_t)+V_{t+1}(F_t(s_t,a_t,\xi_t))].
 \end{equation}
 Because of the assumption that the distribution of the process
$\xi_1,...,\xi_T$, does not depend on the decisions, the set of policies $\Pi$ can be viewed as functions $a_t=\pi_t(\xi_1,...,\xi_{t-1})$ of the data process, $t=1,...,T$, with initial value $\xi_0$. Then for a given (non-randomized)  policy $\pi$, the expectation in \eqref{soc-1} is taken with respect to the probability distribution of $(\xi_1,...,\xi_T)$.
 The  dynamic equations \eqref{contdyn}  correspond to the   problem
 \begin{equation}\label{drcont}
\min\limits_{\pi\in \Pi}   \hat{\cR}\left [ \tsum_{t=1}^{T}
c_t(s_t,a_t,\xi_t) +c_{T+1}(s_{T+1})
\right],\\
 \end{equation}
 where $\Pi$ is the set of the respective non-randomized policies.
Here $\Z=\Z_1\times\cdots\times \Z_T$, where $\Z_t$ is the space of functions $Z_t:\Xi_t\to \bbr$, and
$\hat{\cR}:\Z\to \bbr$ is the nested functional given by the composition of the conditional counterparts of the functionals
\[
 \rho_t(Z_t):=\sup_{Q_t\in \cQ_t}\bbe_{Q_t}[Z_t],\;Z_t\in \Z_t,
 \]
 (cf., \cite{PicSha2021}).

Note   the essential difference between the SOC and MDP approaches to distributionally robust optimization. In the MDP framework,  in the   formulation \eqref{risk-1} the nested functional $\cR^\pi$  is defined on the scenario tree determined by scenarios $s_1,...,s_{T+1}$ of realizations of the state variables, and depends on the (non-randomized) policy $\pi$. On the other hand in the SOC framework the nested functional $\hat{\cR}$ corresponds to the sample paths (scenarios) of the data process $\xi_1,...,\xi_T$, and does not depend on a particular policy $\pi$.

\begin{remark}
{\rm
The ambiguity set can be written here as
\begin{equation}\label{setp-1}
 \PP_t=\left\{P_t: P_t(\cdot|s_t,a_t)=P_t^{Q_t}(\cdot|s_t,a_t),\; Q_t\in \cQ_t,\;  (s_t,a_t)\in \s_t\times \A_t\right\},
\end{equation}
and  the sets    $\cP^t_{s_t,a_t}$ and $\cP^t_{s_t}$, defined in \eqref{sa_marginal} and \eqref{eq_marginal_s}, respectively, can be written   as
\begin{eqnarray}
\label{ambset-1}
\cP^t_{s_t,a_t}=\left\{q_t:q_t(\cdot)=P_t^{Q_t}(\cdot|s_t,a_t)\;\text{for some}\;Q_t\in \cQ_t\right\},\\
\cP^t_{s_t}=\left\{q_t:q_t(\cdot,\cdot)=
 P_t^{Q_t}(\cdot|s_t,\cdot)\;\text{for some}\;Q_t\in \cQ_t\right\},
 \label{ambset-2}
\end{eqnarray}
where $P_t^{Q_t}$ is defined in \eqref{index}.
We have that
   the $(\mathrm{s}, \mathrm{a})$-rectangularity holds iff
    $P_t(\cdot|s_t,a_t)\in  \cP^t_{s_t,a_t}$ for all $P_t\in  \PP_t$ and all
  $ (s_t,a_t)\in \s_t\times \A_t$.
   The $\mathrm{s}$-rectangularity holds iff
    $P_t(\cdot|s_t,\cdot)\in  \cP^t_{s_t}$ for all $P_t\in  \PP_t$ and all
  $ s_t\in \s_t$.
   In particular, suppose that the points $F_t(s_t,a_t,\xi_t)$ are all different from each other  for different $(s_t,a_t,\xi_t)\in \s_t\times\A_t\times \Xi_t$. Then the $\mathrm{s}$-rectangularity (and hence $(\mathrm{s}, \mathrm{a})$-rectangularity) can hold   only if the ambiguity set $\cQ_t$ is the singleton.
}$\hfill \square$
\end{remark}

\section{Conclusions}

For robust MDPs with static formulation of the form \eqref{eq_static_formulation},  the arguably most important factor is its tractability.
Indeed it is shown in \cite{kuhn2013} that for general ambiguity sets the static formulation \eqref{eq_static_formulation} is NP-hard.
Consequently many efforts have been devoted to constructing structured ambiguity sets that allow efficient computation. Most of the time this hinges upon certain dynamic equations or fixed-point characterization of the value function in the stationary case.
The $(\mathrm{s}, \mathrm{a})$-, $\mathrm{s}$-, and $\mathrm{r}$-rectangularity are some notable discoveries among this line of research.
Nevertheless each of them requires a case-by-case analysis in establishing such dynamic equations or the strong duality of the static formulation.
Through our discussions we argue that these properties of static formulation arise precisely as its equivalence to the game formulation and the strong duality of the game formulation.
Such a connection also creates possibilities for new construction of ambiguity sets for static formulation with tractability.
In Section \ref{sec-soc} we show that from the point of view of the distributionally robust setting,   the  Stochastic Optimal Control model  differs from its MDP counterpart in that in seemingly natural formulations the basic condition of  $\mathrm{s}$-rectangularity and $(\mathrm{s}, \mathrm{a})$-rectangularity can hold only if   all ambiguity sets are the singletons.
\\

\noindent{\bf Acknowledgement}\\
The authors are indebted to Eugene   Feinberg for  insightful discussions  which helped to improve the presentation.

\bibliographystyle{plain}
\bibliography{references}

\end{document}